\documentclass [12pt,a4paper,reqno]{amsart}
\usepackage{amsmath,amsthm}
\usepackage{amsfonts}
\usepackage{amssymb}
\allowdisplaybreaks

\newcommand{\be}{\begin{equation}}
\newcommand{\ef}{\end{equation}}
\chardef\bslash=`\\ 

\newtheorem*{thm*}{Theorem}

\theoremstyle{definition}

\newtheorem*{remark*}{Remarks}
\newtheorem*{defn*}{Definition}

\theoremstyle{remark}

\newcommand{\G}{\Gamma}
\newcommand{\wt}{\widetilde}
\newcommand{\wh}{\widehat}
\newcommand{\fc}{\frac}

\newcommand{\iy}{\infty}
 \renewcommand{\sectionmark}[1]{}

\newcommand{\Be}{Beltrami}

\newcommand{\qc} {quasiconformal}

\newcommand{\ve}{\varepsilon}

\newcommand{\const}{\operatorname{const}}

 \usepackage{amsfonts}
\newcommand{\field}[1]{\mathbb{#1}}
\newcommand{\g}{\gamma}

\newcommand{\D}{\field{D}}
\newcommand{\om}{\omega}
\newcommand{\z}{\zeta}
\newcommand{\ov}{\overline}
\newcommand{\vp}{\varphi}
\newcommand{\hC}{\wh{\field{C}}}
\newcommand{\C}{\field{C}}

\newcommand{\B}{\mathbf{B}}
\newcommand{\T}{\mathbf{T}}
\newcommand{\Belt}{\operatorname{Belt}}

\newcommand{\grad}{\operatorname{grad}}
\newcommand{\id}{\operatorname{id}}

\newcommand{\Fib}{\operatorname{Fib}}
\newcommand{\Teich}{\operatorname{Teich}}

\newcommand{\vk} {\varkappa}
\newcommand{\x} {\mathbf x}
\renewcommand{\a} {\alpha}

\newcommand{\ld}{\lambda}
\newcommand{\kp}{\kappa}

\begin{document}

\title{Towards a general distortion theory for univalent functions:
Teichm\"{u}ller spaces and coefficient problems of complex analysis}
\author{Samuel L. Krushkal}

\begin{abstract} Estimating the coefficient functionals on various classes of holomorphic functions
traditionally forms an important field of geometric complex analysis and its mathematical and physical applications.
These coefficients reflect fundamental intrinsic features of holomorphy and of conformality.

This paper surveys the results obtained by a new approach involving deep features of Teichm\"{u}ller spaces. This approach was recently suggested by the author. The paper also contains some new results generalizing the classical coefficient conjectures and presents open problems.
\end{abstract}

\date{\today\hskip4mm({TeSpCoPrCA.tex})}

\maketitle

\bigskip

{\small {\textbf {2020 Mathematics Subject Classification:} Primary: 30C35, 30C50, 30C62, 30C75, 30F60;
Secondary 30C55, 31A05, 32L81, 32Q45, 46G20}

\medskip

\textbf{Key words and phrases:} Teichm\"{u}ller spaces, univalent functions, quasiconformal
extension, coefficient problems, Bieberbach conjecture, nonvanishing holomorphic functions,
subharmonic functions, Schwarzian derivative, holomorphic and plurisubharmonic functionals,
the Hummel-Scheinberg-Zalcman conjecture, Krzyz's conjecture,  Bers' isomorphism theorem,
polynomial approximation, univalent polynomials}

\bigskip

\markboth{Samuel L. Krushkal}{General distortion theory for univalent functions} \pagestyle{headings}

\bigskip\bigskip
\centerline{\bf 1. PREAMBLE}

\bigskip
Estimating the coefficient and more general holomorphic functionals on various classes of holomorphic (especially, of univalent) functions is important in many geometric and physical applications of complex analysis, for example, in view of their connection with string theory and with a holomorphic extension of the Virasoro algebra. The Taylor coefficients reflect key intrinsic features of holomorphic  maps.

This topic has a long history but still remains central in geometric complex function theory and is intensively investigated by many authors using deep methods.

Recently the author presented in \cite{Kr6}, \cite{Kr7}, \cite{Kr15} a new approach to solving the classical coefficient problems on various classes of holomorphic functions, not necessarily univalent, which involves deep analytic and geometric features of Teichm\"{u}ller spaces, especially the Bers isomorphism theorem for Teichm\"{u}ller spaces of punctured Riemann surfaces.

The present paper surveys the results obtained in this way, strengthening some of them, and states open problems.

\newpage
\centerline{\bf 2. TEICHM\"{U}LLER SPACES}

\bigskip
First we recall briefly some needed results from Teichm\"{u}ller space theory on spaces
involved in order to prove our main theorems. A detailed exposition of this theory can be found,
for example, in \cite{Be2}, \cite{GL}, \cite{Kr1}, \cite{L}.
We shall use the notations
$$
\D = \{z \in \C: \ |z| , 1\}, \quad \D^* = \{z \in \hC = \C \cup \{\iy\}: \ |z| > 1\},
\quad \mathbb S^1 = \{z: \ |z| = 1\}.
$$

\bigskip\noindent
{\bf 2.1}. The {\bf universal Teichm\"{u}ller space} $\T = \Teich (\D)$ is the space of
quasisymmetric homeomorphisms of the unit circle $\mathbb S^1$ factorized by
M\"{o}bius maps;  all Teichm\"{u}ller spaces have their isometric copies in $\T$.

The canonical complex Banach structure on $\T$ is defined by factorization of the ball of
the {\bf Beltrami coefficients} (or complex dilatations)
$$
\Belt(\D)_1 = \{\mu \in L_\iy(\C): \ \mu|\D^* = 0, \ \|\mu\| < 1\},
$$
letting $\mu_1, \mu_2 \in \Belt(\D)_1$ be equivalent if the
corresponding \qc \ maps $w^{\mu_1}, w^{\mu_2}$ (solutions to the
Beltrami equation $\partial_{\ov{z}} w = \mu \partial_z w$ with $\mu
= \mu_1, \mu_2$) coincide on the unit circle $\mathbb S^1 = \partial \D^*$
(hence, on $\ov{\D^*}$). Such $\mu$ and the corresponding maps
$w^\mu$ are called $\T$-{\it equivalent}. The equivalence classes
$[w^\mu]_\T$ are in one-to-one correspondence with the {\bf Schwarzian derivatives}
$$
S_w(z) = \left(\frac{w^{\prime\prime}(z)}{w^\prime(z)}\right)^\prime
- \frac{1}{2} \left(\frac{w^{\prime\prime}(z)}{w^\prime(z)}\right)^2
\quad (w = w^\mu(z), \ \ z \in \D^*).
$$

The solutions of the Beltrami equation $\partial_{\ov{z}} w = \mu \partial_z w$ with $\mu \in \Belt(\D)_1$ and of the Schwarz equation $S_w(z) = \vp(z)$ with $\vp \in \B(D^*)$ are conformal
on the disk $\D^*$ and are defined up to a M\"{o}bius (fractional-linear) transformation of $\hC$.
We shall use the hydrodynamical normalization on the infinite point $z = \iy$, i.e., let
$$
w^\mu(z) = z + b_0 + b_1 z^{-1} + b_2 z^{-2} + \dots
$$
and some additional normalization, for example $w^\mu(1) = 1$.
The chain rule
$$
S_{f_1 \circ f}(z) = (S_{f_1} \circ f) f^\prime(z)^2 + S_f(z).
$$
yields for the M\"{o}bius map $w = \g(z)$ the equality
$$
S_{f_1 \circ \g}(z) = (S_{f_1} \circ \g) \g^\prime(z)^2,
\quad S_{\g \circ f}(z) = S_f(z).
$$
Hence, each $S_f(z)$ can be regarded as a quadratic differential $\vp = S_f(z) dz^2$ on $\D$.
T solution $w(z)$ of the Schwarzian equation $S_w(z) = \vp(z)$ with a given holomorphic $\vp$ is defined up to a M\"{o}bius transformation of $\hC$.

Note also that for each locally univalent function $w(z)$ on a simply
connected hyperbolic domain $D \subset \hC$, its Schwarzian
derivative belongs to the complex Banach space $\B(D)$ of
hyperbolically bounded holomorphic functions on $D$ with the norm
$$
\|\vp\|_\B = \sup_D \ld_D^{-2}(z) |\vp(z)|,
$$
where $\ld_D(z) |dz|$ is the hyperbolic metric on $D$ of Gaussian
curvature $- 4$; hence $\vp(z) = O(z^{-4})$ as $z \to \iy$ if $\iy
\in D$. In particular, for the unit disk,
$$
\ld_\D(z) = 1/(1 - |z|^2).
$$

The space $\B(D)$ is dual to the Bergman space $A_1(D)$, a subspace
of $L_1(D)$ formed by integrable holomorphic functions (quadratic differentials
$\vp(z) dz^2$ on $D$), since every linear functional $l(\vp)$ on $A_1(D)$ is
represented in the form
 \be\label{1}
l(\vp) = \langle \psi, \vp \rangle_D = \iint\limits_D \
\ld_D^{-2}(z) \ov{\psi(z)} \vp(z) dx dy
\end{equation}
with a uniquely determined $\psi \in \B(D)$.

The Schwarzians $S_{w^\mu}(z)$ with $\mu \in \Belt(\D)_1$ range over a bounded domain in the space
$\B = \B(\D^*)$. This domain models the space $\T$.
It lies in the ball $\{\|\vp\|_\B < 6\}$ and contains the ball $\{\|\vp\|_\B < 2\}$.
In this model, the Teichm\"{u}ller spaces of all hyperbolic Riemann surfaces are
contained in $\T$ as its complex submanifolds.

The factorizing projection
$$
\phi_\T(\mu) = S_{w^\mu}: \ \Belt(\D)_1 \to \T
$$
is a holomorphic map from $L_\iy(\D)$ to $\B$. This map is a split submersion, which means that
$\phi_\T$ has local holomorphic sections (see, e.g., [GL]).

Note that both equations $S_w = \vp$ and $\partial_{\ov z} w = \mu
\partial_z w$ (on $\D^*$ and $\D$, respectively) determine their solutions in $\Sigma_\theta$ uniquely, so the values $w^\mu(z_0)$ for any fixed $z_0 \in \C$ and the Taylor coefficients $b_1, b_2,
\dots$ of $w^\mu \in \Sigma_\theta$ depend holomorphically on $\mu \in \Belt(\D)_1$ and on
$S_{w^\mu} \in \T$.

\bigskip\noindent
{\bf 2.2}. We shall also use other base points of the space $\T$.
Let $L \subset \C$ be an oriented bounded quasicircle separating the points $0$ and $\iy$. Denote its interior and exterior domains by $D$ and $D^*$ (so $0 \in D, \ \iy \in D^*$).
Then, if $\delta_D(z)$ denotes the Euclidean distance of $z$ from the boundary of $D$ and $\ld_D(z) |dz|$ is its hyperbolic metric of Gaussian curvature $-4$, we have
$$
\fc{1}{4} \le \ld_D(z) \delta_D(z) \le 1.
$$
The right hand inequality follows from the Schwarz lemma and the left from the Koebe one-quoter theorem.
In particular, $\ld_\D(z) = 1/(1 - |z|^2), \ \ld_{\D^*}(z) = 1/(|z|^2 - 1)$.

Consider the classes $\Sigma_{D^*, \theta}$ of univalent functions in the outer domain $D^*$ with expansions
$$
F(z) =  e^{- i \theta} z + b_0 + b_1 z^{-1} + b_2 z^{-2} + \dots \quad (\theta \ \ \text{fixed})
$$
admitting quasiconformal extensions onto $D$. Using $\T$-equivalence of quasiconformal homeomorphisms
$w^\mu$ with $\mu$ from the ball
$$
\Belt(D)_1 = \{\mu \in L_\iy(\C): \ \mu(z)|D^* = 0, \ \ \|\mu\|_\iy  < 1\}.
$$
one obtains in the similar fashion the desired space $\T = \Teich(D)$ with the base point $D^*$,
formed by the Schwarzian derivatives $S_{w^\mu} \in \B(D^*)$.

\bigskip\noindent
{\bf 2.3}. The points of {\bf Teichm\"{u}ller space $\T_1 = \Teich(\D_{*})$ of the punctured disk} $\D_{*} = \D \setminus \{0\}$ are the classes $[\mu]_{\T_1}$ of $\T_1$-{\it equivalent} Beltrami coefficients $\mu \in \Belt(\D)_1$ so that the corresponding quasiconformal  automorphisms $w^\mu$ of the unit disk coincide on both boundary components (unit circle $\mathbb S^1 = \{|z| =1\}$ and the puncture $z = 0$)
and are homotopic on $\D \setminus \{0\}$. This space can be endowed with a canonical complex structure of a complex Banach manifold and embedded into $\T$ using uniformization.

Namely, the disk $\D_{*}$ is conformally equivalent to the factor $\D/\G$, where $\G$ is a cyclic parabolic Fuchsian group acting discontinuously on $\D$ and $\D^*$. The functions $\mu \in L_\iy(\D)$ are lifted to $\D$ as the \Be \ $(-1, 1)$-measurable forms  $\wt \mu d\ov{z}/dz$ in $\D$ with respect to $\G$, i.e., via
$(\wt \mu \circ \g) \ov{\g^\prime}/\g^\prime = \wt \mu, \ \g \in \G$, forming the Banach space $L_\iy(\D, \G)$.

We extend these $\wt \mu$ by zero to $\D^*$ and consider the unit ball $\Belt(\D, \G)_1$ of $L_\iy(\D, \G)$. Then the corresponding Schwarzians $S_{w^{\wt \mu}|\D^*}$ belong to $\T$. Moreover, $\T_1$
is canonically isomorphic to the subspace $\T(\G) = \T \cap \B(\G)$, where $\B(\G)$ consists of elements $\vp \in \B$ satisfying $(\vp \circ \g) (\g^\prime)^2 = \vp$ in $\D^*$ for all $\g \in \G$.

Due to the Bers isomorphism theorem, the space $\T_1$ is biholomorphically isomorphic to the Bers fiber space
$$
\mathcal F(\T) = \{(\phi_\T(\mu), z) \in \T \times \C: \ \mu \in
\Belt(\D)_1, \ z \in w^\mu(\D)\}
$$
over the universal space $\T$ with holomorphic projection $\pi(\psi, z) = \psi$ (see \cite{Be2}).

This fiber space is a bounded hyperbolic domain in $\B \times \C$ and represents the collection of domains $D_\mu = w^\mu(\D)$ as a holomorphic family over the space $\T$. For every $z \in \D$,  its
orbit $w^\mu(z)$ in $\T_1$ is a holomorphic curve over $\T$.

The indicated isomorphism between $\T_1$ and $\mathcal F(\T)$ is induced by the inclusion map \linebreak
$j: \ \D_{*} \hookrightarrow \D$ forgetting the puncture at the origin via
 \be\label{2}
\mu \mapsto (S_{w^{\mu_1}}, w^{\mu_1}(0)) \quad \text{with} \ \
\mu_1 = j_{*} \mu := (\mu \circ j_0) \ov{j_0^\prime}/j_0^\prime,
\end{equation}
where $j_0$ is the lift of $j$ to $\D$.

In the line with our goals, we slightly modified the Bers  construction, applying
quasiconformal maps $F^\mu$ of $\D_{*}$ admitting conformal extension to $\D^*$
(and accordingly using  the Beltrami coefficients $\mu$ supported in the disk)
(cf. \cite{Kr10}). These changes are not essential and do not affect the underlying
features of the Bers isomorphism (giving the same space up to a biholomorphic isomorphism).

The Bers theorem is valid for Teichm\"{u}ller spaces $\T(X_0 \setminus \{x_0\})$ of all punctured hyperbolic Riemann surfaces $X_0 \setminus \{x_0\}$ and implies that $\T(X_0 \setminus \{x_0\})$ is
biholomorphically isomorphic to the Bers fiber space $ \Fib(\T(X_0))$ over $\T(X_0)$.

\bigskip\noindent
{\bf 2.4}. We shall consider also the {\bf space $\T(0, n)$ of punctured spheres} (Riemann surfaces of genus zero)
$$
X_{\mathbf z} = \hC \setminus \{0, 1, z_1 \dots, z_{n-3}, \iy\}
$$
defined by ordered $n$-tuples $\mathbf z = (0, 1, z_1, \dots, z_{n-3}, \iy), \ n > 4$ with distinct $z_j \in \C \setminus \{0, 1\}$.

Fix a collection $\mathbf z^0 = (0, 1, z_1^0, \dots, z_{n-3}^0, \iy)$ with $ z_j^0 \in S^1$ defining the base point $X_{\mathbf z^0}$ of the Teichm\"{u}ller space $\T(0, n) = \T(X_{\mathbf z^0})$. Its
points are the equivalence classes $[\mu]$ of Beltrami coefficients from the ball $\Belt(\C)_1 = \{\mu \in L_\iy(\C): \ \|\mu\|_\iy < 1\}$ under the relation: $\mu_1 \sim \mu_2$, if the corresponding
\qc \ homeomorphisms $w^{\mu_1}, w^{\mu_2}: \ X_{\mathbf a^0} \to X_{\mathbf a}$  are homotopic on $X_{\mathbf a^0}$ (and hence coincide in the points $0, 1, z_1^0, \dots, z_{n-3}^0, \iy$). This
models $\T(0, n)$ as the quotient space $\T(0, n) = \Belt(\C)_1/\sim $ with complex Banach structure of dimension $n - 3$ inherited from the ball $\Belt(\C)_1$.

Another canonical model of $\T(0, n) = \T(X_{\mathbf z^0})$ is obtained again using the uniformization. The surface $X_{\mathbf z^0}$ is conformally equivalent to the quotient space $U/\G_0$, where
$\G_0$ is a torsion free Fuchsian group of the first kind acting discontinuously on $\D \cup \D^*$. The functions $\mu \in L_\iy(X_{\mathbf z^0})$ are lifted to $\D$ as the Beltrami $(-1, 1)$-measurable forms  $\wt \mu d\ov{z}/dz$ in $\D$ with respect to $\G_0$ which satisfy $(\wt \mu \circ \g) \ov{\g^\prime}/\g^\prime = \wt \mu, \ \g \in \G_0$ and form the Banach space $L_\iy(\D, \G_0)$.

After extending these $\wt \mu$ by zero to $\D^*$, the  Schwarzians $S_{w^{\wt \mu}|\D^*}$ for $\|\wt \mu\|_\iy$  belong to $\T$ and form its subspace regarded as the {\it Teichm\"{u}ller space
$\T(\G_0)$ of the group $\G_0$}. It is canonically isomorphic to the space $\T(X_{\mathbf z^0})$, and moreover,
$$
\T(\G_0) = \T \cap \B(\G_0),
$$
where $\B(\G_0)$ is an $(n - 3)$-dimensional subspace of $\B$ which consists of elements $\vp \in \B$ satisfying $(\vp \circ \g) (\g^\prime)^2 = \vp$ for all $\g \in \G_0$ (holomorphic
$\G_0$-automorphic forms of degree $- 4$); see, e.g. \cite{Le}.

This leads to the representation of the space $\T(X_{\mathbf z^0})$ as a bounded domain in the complex Euclidean space $\C^{n-3}$.

Note that $\B(\G_0)$ has the same elements as the space $A_1(\D^*, \G_0)$ of integrable holomorphic forms of degree $- 4$ with norm $\|\vp\|_{A_1(\D^*, \G_0)} = \iint_{\D^*/\G_0} |\vp(z)| dx dy$; and
similar to (1), every linear functional $l(\vp)$ on $A_1(\D^*, \G_0)$ is represented in the form
$$ l(\vp) = \langle \psi, \vp \rangle_{\D/\G_0} :=
\iint\limits_{\D^*/\G_0} (1 - |z|^2)^2 \ \ov{\psi(z)} \vp(z) dx dy
$$
with uniquely determined $\psi \in \B(\G_0)$.

Any Teichm\"{u}ller space is a complete metric space with intrinsic Teichm\"{u}ller metric defined by quasiconformal maps. By the Royden-Gardiner theorem, this metric equals the hyperbolic  Kobayashi metric
determined by the complex structure (see, e.g., \cite{EKK}, \cite{GL}, \cite{Ro}).

\bigskip\noindent
{\bf 2.5. Weak approximation of Teichm\"{u}ller spaces}.
We apply a weak (locally uniform) approximation of the underlying space $\T$ and simultaneously
of the space $\T_1$ by finite dimensional Teichm\"{u}ller spaces of the punctured spheres arised after
 exhaustion of the unit disk  by the canonical fundamental polygons of these surfaces and from the locally uniform convergence of associated univalent functions on $\C$.

Note that the space $\B$ containing $\T$ is not separable;  thus, it does not admit a strong sequential approximation in $\B$-norm. However, the weak (locally uniform) convergence of the Schwarzians $S_f$
on $\D^*$ (together with $f \in \Sigma$) also restrains the growth of the lifted functional $J$ to $\Fib  \T$ and ensures that the maximum of $|\mathcal J|$ cannot increase under applied semicontinuous regularization in the strong topology on $\T$.

\bigskip\bigskip
\centerline{\bf 3. SHARP COEFFICIENT ESTIMATES FOR ROTATIONALLY INVARIANT}
\centerline{\bf COLLECTIONS OF UNIVALENT FUNCTIONS}

\bigskip\noindent
{\bf 3.1. Classes of functions}. Consider the collection of univalent functions
$$
w(z) = a_1 z + a_2 z^2 + \dots \quad \text{with} \ \ |a_1| = 1,
$$
on the unit disk $\D$ admitting quasiconformal extensions to the whole Riemann sphere $\hC$ and assume
that these extensions satisfy $w(1) = 1$. Take the completion of this collection in the topology of
locally uniform convergence on $\D$ and denote it by $\wh S(1)$. It is a disjunct union
$$
\wh S(1) = \bigcup_{- \pi \le \theta < \pi} S_\theta(1),
$$
of sets $S_\theta(1)$ of univalent functions $w(z) = e^{i \theta} z + a_2 z^2 + \dots$ with fixed $\theta$. In the general case, the equality $w(1) = 1$ must be understand (for the limit functions
of sequences of the initial quasiconformally extendible functions) in terms of the Carath\'{e}odory
prime ends.

The family $\wh S(1)$ closely relates to the canonical class $S$ of univalent functions $w(z$
on the unit disk normalized by $w(0) = 0, \ w^\prime(0) = 1$. Every $w \in S$ has its representative $\wh w$ in $\wh S(1)$ (not necessarily unique) obtained by pre and post compositions of $w$ with rotations $z \mapsto e^{i \a} z$ about the origin, related by
 \be\label{3}
w_{\a, \beta}(z) =  e^{- i \beta} w(e^{i \a} z) \quad \text{with} \ \ \a = \arg z_0,
\end{equation}
where $z_0$ is a point whose image $w(z_0) = e^{i \theta}$ is a common
point of the unit circle and the  boundary of domain $w(\D)$. The existence of such point follows
from Schwarz's lemma.

This implies, in particular, that the functions conformal in the closed disk
$\ov \D$ are dense in each class $S_\theta(1)$. Such a dense subset is formed, for example, by the
images of the homotopy functions $[f]_r(z) = \fc{1}{r} f(r z)$ with real $r \in (0, 1)$ combined
with rotations (1).

Our goal is to estimate sharply on rather general subclasses $\mathcal X$ of $S$ the general polynomial coefficient functionals
\be\label{4}
J(f) = J(a_{m_1}, \dots, a_{m_s}), \quad 2 < a_{m_1} < \dots < a_{m_s} < \iy,
\end{equation}
depending from a distinguished finite set of coefficients $a_{m_j}(f)$.
We assume that these functionals and the indicated subclasses $\mathcal X$ are invariant under rotations (3) with independent $\a$ and $\beta$ from $[- \pi, \pi]$); so $|J(f_{\a, \beta})| = |J(f)|$).

The inverted functions $F_f(z) = 1/f(1/z)$ belong to the class $\Sigma$ of univalent $\hC$-holomorphic functions on the disk $\D^*$ with expansions
$$
F(z) = z + b_0 + b_1 z^{-1} + b_2 z^{-2} + \dots
$$
(such functions have been applied above in construction of Teichm\"{u}ller space $\T$).
Denote $S_Q$ and $\Sigma_Q$ the (dense) subclasses of $S$ and $\Sigma$ formed by functions having quasiconformal extension to $\hC$.

Now pick  a rotationally invariant subclass $\mathcal X$ of $\wh S(1)$, which satisfies the following  conditions:

$(a)$ {\bf openness}, which means that the corresponding collection of the Schwarzians $S_w$  of $w \in \mathcal X$ determines in the space $\T$  a complex Banach
submanifold (of finite or infinite dimension), which we denote by $\mathcal X_\T$;

$(b)$ {\bf variational stability}, which means that for any quasiconformal deformation $h$  of a functions $w \in \mathcal X$ with Beltrami coefficients $\mu_h$ (and sufficienly small dilatation $\|\mu\|_\iy$)   supported in the complementary domain of $w(\D)$ the composition $h \circ w|\D$ also belongs to $\mathcal X$.

In fact, we shall use only a special type of such quasiconformal deformations.

We associate with such a class $\mathcal X$ the quantity
  \be\label{5}
|a_2(\mathcal X)| = \max \{|a_2(w)|: \ w \in \mathcal X\}
\end{equation}
and the set of rotations
 \be\label{6}
\mathcal R_\mathcal X  = \{w_{0,\tau,\theta}(z) = e^{- i \theta} w_0(e^{i \tau} z)\},
\end{equation}
where $w_0$ is one of the maximizing functions for $a_2$
on $\mathcal X$, i.e., with $|a_2(w_0)| = |a_2(\mathcal X)|$.

We prove for such classes of univalent functions and functionals a general theorem, which completely describes the extremal these functionals. This theorem covers the known general distortion results
for coefficients.

\bigskip\noindent
{\bf 3.2. Underlying theorem}. We precede the main theorem by the following underlying result, which
reveals the intrinsic connection between the covering and distortion features of holomorphic functions.

Let $G$ be a domain in a complex Banach space $X = \{\mathbf x\}$
and $\chi$ be a holomorphic map from $G$ into the universal Teichm\"{u}ller space $\T = \Teich(D)$ with the base point $D$ modeled as a bounded subdomain of $\B$. Assume that $\chi(G)$ is a (pathwise
connected) submanifold of finite or infinite dimension in $\T(D)$, and put, in accordance with (4),
 \be\label{7}
a_{2,\theta}^0 = |a_2(\mathcal \chi(G))| = \sup \{|a_2|: \ S_w \in \chi(G)\},
\end{equation}

\noindent
{\bf Theorem 1}. \cite{Kr15} {\it Let $w(z)$ be a holomorphic univalent solution of the Schwarz differential equation
$$
S_w(z) = \chi(\x)
$$
on $D$ satisfying $w(0) = 0, \ w^\prime(0) = e^{i \theta}$ with the fixed $\theta \in [-\pi, \pi]$
and $\x \in G$ (hence $w(z) = e^{i \theta} z  + \sum_2^\infty a_n z^n$). Assume that
$a_{2,\theta}^0 \ne 0$, and let $w_0(z) = e^{i \theta} z + a_2^0 z^2 + \dots$ be one of the maximizing functions for $a_{2,\theta}^0$. Then:

(a) For every indicated function $w(z)$ , the image domain $w(D)$ covers entirely the disk
$D_{1/(2 |a_{2,\theta}^0|)} = \{|w| < 1/(2 |a_{2,\theta}^0|)\}$.

The radius value $1/(2 |a_{2, \theta}^0|)$ is sharp for this collection of
functions and fixed $\theta$, and the circle $\{|w| = 1/(2 |a_{2,\theta}^0|)$ contains points
not belonging to $w(\D)$ if and only if $|a_2| = |a_{2,\theta}^0|$
(i.e., when $w$ is one of the maximizing functions).

(b) The inverted functions
$$
W(\zeta) = 1/w(1/\zeta) = e^{i \theta}\zeta - a_2^0 + b_1 \zeta^{-1} + b_2 \zeta^{-2} + \dots
$$
with $\z \in D^{-1}$ map domain $D^{-1}$ onto a domain whose boundary is entirely
contained in the disk} $\{|W + a_{2,\theta}^0| \le |a_{2,\theta}^0|\}$.

{\it (c) Let $D$ be the unit disk $\D$ and
$$
\vp(z) = \chi(\x) = c_0 + c_1 z + c_2 z^2 + \dots, \quad |z| < 1.
$$
Let a point $\x_0 \in G$ correspond to the function $w_0(z)$ maximizing the quantity $|a_{2,\theta}^0|$ on $[-\pi, \pi]$. If its Schwarzian
$S_{w_0} = \chi(\x_0) = c_0^0 + c_1^0 + \dots$ satisfies
$$
c_1^0 \ne 0,
$$
then the value $|c_1^0|$ is maximal on this class, i.e.,
$$
|c_1^0| = \sup \ \{|c_1(w_{\a,\beta})|: \ S_w \in \chi(G); \ \a, \beta \in [- \pi, \pi]\}.
$$
}

\noindent
{\bf Proof}. The proof of parts {\it (a), (b)} follows the classical lines of Koebe's $1/4$
theorem (cf. \cite{Go}).

The proof of assertion (c) is complicated and actually represents a special case of results stated it
the next theorems.

\bigskip\noindent
{\bf 3.3. Distortion theorems for holomorphic functionals}.
To have the rotational symmetry (3), one must take $L = \mathbb S^1$ and hence, $D = \D,
\ D^* = \D^*$, and deal with the corresponding classes $S_\theta(1)$ (or $\Sigma_\theta(1)$ of inverted functions $F_f$on $\D^*$) and their appropriate subclasses.

The assumption $(b)$ on images $\chi(G)$ provides that $\chi(G)$ is a subdomain of $\T$ having the common boundary points with $\partial \T$ and that the maximal value of $|a_2|$ on $\chi(G)$ is attained on $\partial \T$.

\bigskip\noindent
{\bf Theorem 2}. {\it Any rotationally invariant polynomial functional (4),
whose zero set $\mathcal Z_J = \{w \in \wh S: \ J(w) = 0\}$
is separated from the rotation set (6), is maximized on the class $\mathcal X$ only by  functions $w_{0,\tau,\theta} \in \mathcal R_\mathcal X$.

In other words, any extremal function $w_0$ of any homogeneous (rotationally invariant) coefficient functional $J$ on a rotationally invariant and variationally stable
class $\mathcal X$ must be simultaneously maximal for the second coefficient $a_2$ on this class, unless $J(w_0) = 0$. }

\bigskip
All assumptions of this theorem on the class $\mathcal X$ and the functional $J$
are essential and cannot be omitted.
This will be illustrated on examples below.

In particular, all this holds for the functionals
$$
J(f) = a_n(f), \quad n \ge 3.
$$
on the class $S^0$. Then Theorem 2 implies

\bigskip\noindent
{\bf Theorem 3}. {\it For all functions $f(z) \in \ov{S^0}$ and all $n \ge 3$, we have the sharp estimates
$$
\max_{\ov{S^0}} |a_n| = |a_n(w_0)|,
$$
with equality for each $n$ only for the function $w_0(z)$ maximizing $|a_2|$ on $\ov{S^0}$ .   }

\bigskip
This and the next theorems strengthens the de Branges theorem proving the Bieberbach conjecture \cite{DB}.

\bigskip\noindent
{\bf Theorem 4}. {\it If $\ov{\mathcal S_Q} = S$, then the extremal function for any functional
$J(f)$  of type (3) is the Koebe function
$$
\kp_\theta(z) = \fc{z}{(1 - e^{i \theta} z)^2} = z +
\sum\limits_2^\iy n e^{- i(n-1) \theta} z^n, \quad - \pi < \theta \le \pi.
$$
}

Recall that $\kp_\theta$ maps the unit disk onto the complement of the ray
$$
w = -t e^{-i \theta}, \ \ \fc{1}{4} \le t \le \iy.
$$

\noindent
{{\bf 3.4. Two illustrating examples}.

\noindent
{\bf Example 1}. Consider the balls $B_r = \{\|\vp\|_2 < r\}$ in the Hilbert space $A_2(\D)$ of the
square integrable holomorphic functions
$\vp(z) = c_0 + c_1 z + c_2 z^2 + \dots, \ |z| < 1$,
with norm
$\|\vp\|_2 = \Bigl(\fc{1}{2\pi} \iint\limits_\D |\om|^2 dx dy\Bigr)^{1/2}$.
The Taylor coefficients $c_n$ of these functions coincide with their Fourier coefficients with respect
to the orthonormal system
$$
\vp_n(z) = \sqrt{(n + 1)/\pi} \ z ^ n, \quad n = 0, 1, 2, \dots , \quad \text{on} \ \ \D,
$$
and $\sum\limits_0^\iy |c_n|^2 = \|\vp\|_2$.

In this case, using quasiconformal deformations preserving $L_2$ norm (given, for example, in \cite{Kr7}), one can omit the assumption of Theorem 2 that domain $\chi(G)$ must touch the boundary
of $\T$.

Applying the Schwarz inequality, one obtains
$$
\fc{1}{\pi} \iint\limits_\D |\vp| dx dy \le \fc{1}{\sqrt{\pi}}
\Bigl(\iint\limits_\D |\vp|^2 dx dy\Bigr)^{1/2} = \sqrt{2}   \pi \ \|\vp\|_2,
$$
which together with the inequality
 \be\label{8}
\|\vp\|_\B \le \iint\limits_{\D} |\vp| dx dy
\end{equation}
(following from the mean value inequality for integrable holomorphic functions in a domain) implies
that all functions $\vp \in A_2(\D)$ belong to the space $\B$, and their $\B$-norm is estimated by
 \be\label{9}
\|\vp\|_\B \le \sqrt{2}\pi \  \|\vp\|_2.
\end{equation}
Hence, these functions (more precisely, the corresponding quadratic differentials $\vp(z) dz^2$) can be regarded as the Schwarzian derivatives of locally univalent functions in $\D$.

Now, take the domain $G = B_{r_0}$ with
$$
r_0 = \sqrt{2}/\pi.
$$
By (9), all $\vp \in G$ have the $\B$-norm less than $2$, and by the Ahlfors-Weill theorem \cite{AW}
are the Schwarzian derivatives of univalent functions $w(z)$ on the disk $\D$ having canonical quasiconformal extensions to $\hC$.
Hence, every coefficient $c_n$ is represented as a polynomial from the initial coefficients
$a_2, \dots, a_{n+1}$ of $w(z)$.
The maximal value of each $|c_n|$ equals $1$ and is attained on $G$ on the function $\vp_n$.

Theorems 1 and 2 yield that the univalent function $w_0(z)$ with $S_{w_0}= \vp_1$ maximizes simultaneously both coefficients $a_2$ and $c_1$ on this class, and this function also is extremal
for all $a_n, \ n \ge 3$.

To determine explicitly the Taylor coefficients $a_n$ of univalent $w(z)$ whose Schwarzians run over the closed ball $\ov{B_{r_0}}$, one has to find the ratios $w = \eta_1/\eta_2$ of two independent solutions of the linear differential equation
$$
2 \eta^{\prime\prime}(z) + \vp(z) \eta(z) = 0
$$
with a given holomorphic $\vp$ (which is equivalent to solving the Schwarzian equation $S_w(z) = \vp(z)$)
subject to $w(0) = 0, \ w^\prime(0) =1 $.
The arising linear equations are of the form
$$
\eta^{\prime\prime} + c z^p \ \eta = 0
$$
and closely relate to special functions (see, e.g., \cite{Ka}). For the extremal function $w_0$
distinguished above, we have the equation
$$
\eta^{\prime\prime} + r_0 z \ \eta = 0;
$$
its solutions are represented by linear combinations of cylindrical functions.

\bigskip
Note also that, due to the existence theorem from \cite{Kr3}, for every given $n \in \mathbb N$ and
every bounded holomorphic function $\vp(z) = c_0 + c_1 z + \dots \in A_2(\D)$, which is not a polynomial of degree $n_1 \le n$, there exists a quasiconformal deformation of the extended plane $\hC$ conformal
on $\vp(\D)$, which perturbs arbitrarily (in restricted limits) the coefficients
$c_0, c_1, \dots, c_{n-1}$ but preserves the $L_2$-norm of $\vp$.
The maximizing functions (15) do not admit such deformations.

\bigskip\noindent
{\bf Example 2}.
The class $S(M)$ formed by bounded functions $f(z) = z + a_2 z^2 + \dots \in S$
with $|f(z)| < M$ in $\D \ (M > 1)$ also is not variationally stable, because
the variations given by Lemma 5 (with the sets $E$ of quasiconformality located
outside of the disk $\D_M = \{|z| < M|\}$) generically increase the $\sup$ norm of varied functions.

Thus Theorem 2 cannot be applied to this class. This is in accordance with the fact that the known coefficient estimates for $S(M)$ (see, e.g. \cite{Pr}) are of completely different nature than Theorem 2.

\bigskip
The same reasons cause that this theorem does not work for functions with $k$-quasiconformal extensions
with fixed $k < 1$.

\bigskip\bigskip
\centerline{\bf 4. SKETCH OF PROOF OF THEOREM 2}

\bigskip
The proof of this main theorem is divided to several steps and is based on two lemmas.

\bigskip\noindent
$\mathbf{1^0}$. The following important lemma from \cite{Kr12} provides the existence of quasiconformally extendable functions with some non-standard  normalizations,  which also ensures their  compactness, holomorphic dependence on complex parameters, etc.

\bigskip\noindent
{\bf Lemma 1}. {\it For any Beltrami coefficient $\mu \in \Belt(\D^*)_1$ and any $\theta_0 \in [0, 2 \pi]$, there exists a point $z_0 = e^{i \a}$ located on $\mathbb S^1$ so that
$|e^{i \theta_0} - e^{i \a}| < 1$ and such that for any $\theta$ satisfying
$|e^{i \theta} - e^{i \a}| < 1$ the equation
$\partial_{\ov z} w =  \mu(z) \partial_z w$
has a unique homeomorphic solution $w = w^\mu(z)$, which is holomorphic on the unit disk $\D$
and satisfies
 \be\label{10}
w(0) = 0, \quad w^\prime(0) = e^{i \theta}, \quad w(z_0) = z_0.
\end{equation}
Hence, $w^\mu(z)$ is conformal and does not have a pole in $\D$ \ (so
$w^\mu(z_{*}) = \iy$ at some point $z_{*}$ with $|z_{*}| \ge 1$).  }

\bigskip
The initial point $\mu(z) = 0$ (almost everywhere in $\D^*$) corresponds to the elliptic M\"{o}bius map
$$
\frac{w - z_0}{w} = e^{- i\theta} \ \frac{z - z_0}{z},
$$

This lemma allows one to define the Teichm\"{u}ller spaces using, in particular, the quasiconformally extendible  univalent functions $w(z)$ on $\D$ normalized by
$w(0) = 0, \ w^\prime(0) = 1, \ w(1) = 1$,
and with more general normalization
$$
w(0) = 0, \quad w^\prime(0) = e^{i \theta}, \quad w(1) = 1
$$
(with appropriate $\theta$).
All such functions are holomorphic (have no poles) in the disk $\D$.

Now consider classes $S_\theta(1)$ and their union $\wh S(1)$ introduced in \textbf{3.1} and pass
to the inverted functions $F_f(z) = 1/f(1/z)$ for $f \in \wh S^0$, which form the corresponding classes $\Sigma_{z_0,\theta}$ of nonvanishing  univalent functions on the disk $\D^*$ with expansions
$$
F(z) =  e^{- i \theta} z + b_0 + b_1 z^{-1} + b_2 z^{-2} + \dots,
\quad  F(z_0) = z_0,
$$
and
$$
\Sigma^0 = \bigcup_{z_0,\theta} \Sigma_{z_0,\theta}.
$$
The coefficients $a_n$ of $f(z)$ and the corresponding coefficients $b_j$ of $F_f(z)$ are related by
$$
b_0 + e^{2i \theta} a_2 = 0, \quad b_n + \sum \limits_{j=1}^{n}
\epsilon_{n,j}  b_{n-j} a_{j+1} + \epsilon_{n+2,0} a_{n+2} = 0,
\quad n = 1, 2, ... \ ,
$$
where $\epsilon_{n,j}$ are the entire powers of $e^{i \theta}$. This
successively implies the representations of $a_n$ by $b_j$ via
 \be\label{11}
a_n = (- 1)^{n-1} \epsilon_{n-1,0}  b_0^{n-1} - (- 1)^{n-1} (n - 2)
\epsilon_{1,n-3} b_1 b_0^{n-3} + \text{lower terms with respect to} \ b_0.
\end{equation}
This transforms the initial functional (4) into a coefficient functional $\wt J(F^\mu)$
on $\Sigma^0$ depending on the corresponding coefficients $b_j$. This dependence is holomorphic
from the Beltrami coefficients $\mu_W \in \Belt(\D)_1$ and from the Schwarzians $S_{W^\mu}$.

Accordingly, we model the universal Teichm\"{u}ller space $\T$ by a domain in the space $\B = \B(\D^*)$ formed by the Schwarzians $S_F^\mu$. Thereby, both functionals  $\wt J(F)$ and $J(f)$ are lifted holomorphically onto the space $\T$.

To lift $J$ onto the covering space $\T_1$, we again pass to functional $\wh J(\mu) = \wt J(F^\mu)$
lifting $J$ onto the ball $\Belt(\D)_1$ and apply the $\T_1$-equivalence, i.e., the quotient map
$$
\phi_{\T_1}: \ \Belt(\D)_1 \to \T_1, \quad \mu \to [\mu]_{\T_1}.
$$
Thereby the functional $\wt J(F^\mu)$ is pushed down to a bounded holomorphic functional $\mathcal J$
on the space $\T_1$ with the same range domain.

Using the Bers isomorphism theorem, one can regard the points of the space $\T_1$ as the pairs
$X_{F^\mu} = (S_{F^\mu}, F^\mu(0))$, where $\mu \in \Belt(\D)_1$ obey $\T_1$-equivalence (hence, also $\T$-equivalence). Note that since the coefficients $b_0, \ b_1, \dots$ of $F^\mu \in \Sigma_\theta$   are uniquely determined by its Schwarzian $S_{F^\mu}$, the values of $\mathcal J$ in the
points $X_1, \ X_2 \in \T_1$ with $\iota_1(X_1) = \iota_1(X_2)$ are equal.

In result, we get on $\T_1 = \Fib (\T)$ the holomorphic functional
$$
\mathcal J(X_{F^\mu}) = \mathcal J(S_{F^\mu}, \ t), \quad t = F^\mu(0),
$$
and must investigate the restriction of plurisubharmonic functional $|\mathcal J(S_{F^\mu}, t)|$ to
the image of domain $\chi(G)$ in $\Fib (\T)$. .

Note that by the part {\it (b)} of Theorem 1, the boundary of domain $W^\mu(\D^*)$ under any function
$W^\mu(z) \in \Sigma^0$ is located in the disk $\{|W - b_0| \le |a_2^0| \}$  with $|a_2^0|$
determined by (3). For all such $W^\mu$, the variable $t$ in the representation (11) runs over some
subdomain $D_1$ in the disk $\D_4 = \{|t| < 4\}$ containing the origin (this subdomain depends on $z_1$).
Since the functional $J$ is rotationally invariant, this subdomain $D_1$ is a disk $\D_{\a_1}$ of some radius $\a_1 \le 2 |a_2^0|$.

We define on this domain the function
 \be\label{12}
\wt u_1(t) = \sup_{S_{W^\mu}} \mathcal J(S_{W^\mu}, t).
\end{equation}
taking the supremum over all $S_{W^\mu} \in \T$ admissible for a given $t = W^\mu(0) \in D_{\a_n}$, that means over the pairs $(S_{W^\mu}, t) \in \Fib(\T)$ with a fixed $t$ and pass to the upper semicontinuous
regularization
$$
u_1(t) = \limsup\limits_{t^\prime \to t} \wt u_1(t^\prime).
$$

\noindent
$\mathbf{2^0}$. Now the crucial step in the proof of Theorem 2 is to establish that the function (12) inherits subharmonicity and find the maximal domain for $t = W^\mu(0)$, where it is defined.
Select on the unit circle a dense subset
$$
\mathbf e = \{z_1, z_2, \dots, z_n, \dots\}, \quad z_1 = e^{i \theta_1}.
$$
Repeating successively for the above construction with fix points $z_1, z_2, \dots$, one obtains  similar
to (12) the corresponding functions $u_1(t), u_2(t), \dots$. Take their upper envelope
$\sup_n u_n(t)$ and its upper semicontinuous regularization $u(t)$.

The second basic lemma states:

\bigskip\noindent
{\bf Lemma 2}. {\it The functions $u(t)$ is logarithmically subharmonic in some domain $D_0 \subset \D$. }

\bigskip
The proof of this lemma given in \cite{Kr9} involves a weak approximation of the underlying
space $\T$ and simultaneously of the space $\T_1$ by finite dimensional Teichm\"{u}ller spaces of
the punctured spheres (in the topology of locally uniform convergence on $\C$, as was described
in {\bf 2.5}).

\noindent
$\mathbf{3^0}$. To find the maximal domain $D_0$ (the range domain of values $t = W^\mu(0)$) and describe its boundary points, we apply,  similar to the above classes $S_\theta(1)$, the classes $S_{z_0,\theta}$ consisting of $f \in S_{Q,\theta}$
with a fix point at $z_0 \in \mathbb S^1$ for all $z_0$ and take their unions
$$
\bigcup_{z_0 \in \mathbb S^1, \theta \in [-\pi,pi]} S_{z_0,\theta}.
$$

Now consider the increasing unions of the quotient spaces
\be\label{13}
\mathcal T_s = \bigcup_{j=1}^s \ \wh \Sigma_{\theta_j}^0/\thicksim \
= \bigcup_{j=1}^s \{(S_{W_{\theta_j}}, W_\theta^\mu(0)) \} \ \simeq
\T_1 \cup \dots \cup \T_1,
\end{equation}
where $\theta_j$ run over a dense subset $\Theta \subset [-\pi, \pi]$, the equivalence relation
$\thicksim$ means $\T_1$-equivalence on a dense subset $\wh \Sigma^0(1)$ in the union $\wh \Sigma(1)$
formed by univalent functions $W_{\theta_j}(z) = e^{-i \theta_j} z + b_0 + b_1 z^{-2} + \dots$ on
$\D^*$ with quasiconformal extension to $\hC$ satisfying $W_{\theta_j}(1) = 1$, and
$$
\mathbf W_\theta^\mu(0) := (W_{\theta_1}^{\mu_1}(0), \dots ,
W_{\theta_s}^{\mu_s}(0)).
$$
The Beltrami coefficients  $\mu_j \in \Belt(\D)_1$ are chosen here independently. The corresponding collection
$$
\beta = (\beta_1, \dots, \beta_s)
$$
of the Bers isomorphisms
$$
\beta_j: \ \{(S_{W_{\theta_j}}, W_{\theta_j}^{\mu_j}(0))\} \to
\mathcal F(\T)
$$
determines a holomorphic surjection of the space $\mathcal T_s$ onto $\mathcal F(\T)$. The maximal
function (12) is determined by
$$
u(t) = \sup_{\theta} u_{\theta_s}(t),
$$
where $u_{\theta_s}$ is obtained by maximization of type (12) over $\mathcal T_s$.

In view of the rotational symmetry of the initial functional $J(f)$ and of the union all spaces (13),
the above construction implies a maximal subharmonic function $u(t)$ on a disk $\D_a$ with $a \le 4$.
(cf. \cite{Kr9}).

Now we use the assumption of Theorem 2 that $\chi(G)$ is a domain in $\T$. Together with the features
of construction of the Bers fiber space $\Fib(\T)$, this implies that image of $\chi(G)$ in the space $\Fib(\T)$ also is a connected submanifold covering $\chi(G)$. Hence, the restriction of the function $u(t)$ constructed above to this image is subharmonic in some domain $D_0^\prime$ depending on
$\chi(G)$.

\bigskip\noindent
$\mathbf{4^0}$. Now we use the assumption of Theorem 2 that $\chi(G)$ is a domain in $\T$. Together with the features of the Bers fiber space $\Fib(\T)$, this implies that image of $\chi(G)$ in the space $\Fib(\T)$ also is a connected submanifold covering $\chi(G)$. Hence, the restriction of the function $u(t)$ constructed above to this image is subharmonic in some domain $D_0^\prime$ depending on
$\chi(G)$.

To find explicitly this domain and its boundary point, we use the second assumption on $\chi(G)$, namely,
its variational stability and apply the following variational lemma, which is a special case of the
general deformations constructed in \cite{Kr1}.

\bigskip\noindent
{\bf Lemma 3}. {\it Let $D$ be a simply connected domain on the Riemann sphere
$\hC$. Assume that there are a set $E$ of positive two-dimensional
Lebesgue measure and a finite number of points
 $z_1, z_2, ..., z_m$ distinguished in $D$. Let
$\a_1, \a_2, ..., \a_m$ be non-negative integers assigned to $z_1,
z_2, ..., z_m$, respectively, so that $\a_j = 0$ if $z_j \in E$.

Then, for a sufficiently small $\ve_0 > 0$ and $\varepsilon \in (0,
\varepsilon_0)$, and for any given collection of numbers $w_{sj}, s
= 0, 1, ..., \a_j, \ j = 1,2, ..., m$ which satisfy the conditions
$w_{0j} \in D$, \
$$
|w_{0j} - z_j| \le \ve, \ \ |w_{1j} - 1| \le \ve, \ \ |w_{sj}| \le
\ve \ (s = 0, 1, \dots   a_j, \ j = 1, ..., m),
$$
there exists a quasiconformal automorphism $h$ of $D$ which is conformal on $D
\setminus E$ and satisfies
$$
h^{(s)}(z_j) = w_{sj} \quad \text{for all} \ s =0, 1, ..., \a_j, \ j
= 1, ..., m.
$$
Moreover, the Beltrami coefficient $\mu_h(z) = \partial_{\bar z}
h/\partial_z h$ of $h$ on $E$ satisfies $\| \mu_h \|_\iy \leq M
\ve$. The constants $\ve_0$ and $M$ depend only upon the sets $D, E$
and the vectors $(z_1, ..., z_m)$ and $(\a_1, ..., \a_m)$.

If the boundary $\partial D$ is Jordan or is $C^{l + \a}$-smooth,
where $0 < \a < 1$ and $l \geq 1$, we can also take $z_j \in
\partial D$ with $\a_j = 0$ or $\a_j \leq l$, respectively.   }

This lemma and Theorem 1 imply that the indicated domain $D_0^\prime$ is rotationally symmetric and connected, hence a disk whose radius must be equal to $2|a_2^0|$.
The maximal function $u(t) = \max\limits_{\mathcal X} |J(f^\mu)|$ is circularly symmetric, hence,
$u(t) = u(|t|)$ and must attains its maximal value on the circle $\{|t| = |2|a_2^0|$.
This completes the proof of the theorem.

\bigskip\bigskip
\centerline{\bf 5. EXTREMALS OF SUBHARMONIC FUNCTIONALS}

\bigskip
The results obtained above for holomorphic functionals can be extended to more general, in particular,
plurisubharmonic functionals. We consider here two types of such functionals.
Other types of subharmonic functionals will be applied in the next section.

\bigskip\noindent
{\bf 5.1. Strengthening Theorem 2}. First we show how the previous results are extended
to a broad class of plurisubharmonic coefficient functionals, which allows one to omit  the rigid constrains to holomorphic functionals such as their rotational invariance and restrictions to location
of the zero set, etc. As a consequence, this extension provides the new extremal features of the Koebe function.

Consider on the class $S$  the general plurisubharmonic polynomial coefficient functionals
of the form
 \be\label{14}
 J(f) = J(\a_{n_1,\dots,n_m}) = \max_{-\pi \le \vartheta_1, \dots, \vartheta_m \le \pi} \ \Big\vert \sum\limits_{|n| = 3}^N \ \a_{n_1,\dots,n_m} a_{n_1} \dots a_{n_m} \  e^{i(n_1 \vartheta_1 + \dots + n_m \vartheta_m)} \Big\vert,
 \end{equation}
where $|n| = n_1 + \dots + n_m$ and the factors $\a_{n_1,\dots,n_m}$ do not depend on
coefficients $a_j$.
For such functionals, we have the following general distortion theorem.

\bigskip\noindent
{\bf Theorem 5}. {\it Any nonconstant plurisubharmonic coefficient functional (14) is maximized on the class $S$ only by the rotated Koebe function
 \be\label{15}
\kp_{\tau, \theta}(z) =  e^{- i \theta} \kp_0(e^{i \tau} z).
\end{equation}
}

In particular, this holds for functionals
$$
J_N(f) = J_N(\a_{n_1,\dots,n_m}) = \max_{-\pi \le \vartheta_1, \dots, \vartheta_m \le \pi} \ \Big\vert \sum\limits_{n_1 + \dots + n_m = N} \ \a_{n_1,\dots,n_m} a_{n_1} \dots a_{n_m} \  e^{i(n_1 \vartheta_1 + \dots + n_m \vartheta_m)} \Big\vert
$$
straightforwardly generalizing $J(f) = a_n$.

The proof of this theorem is given in the same fashion as of Theorem 2.

Theorem 5 also can be strengthened to an assertion of type Theorems 1-2, i.e., for smaller subcollections of functions.

As an interesting consequence, we have

\bigskip\noindent
{\bf Corollary 1}. {\it The Koebe function (and its rotations (15)) maximizes every
trigonometric polynomial
$$
 \sum\limits_{|n| = 3}^N \ a_{n_1} \dots a_{n_m} \ e^{i(n_1 \vartheta_1 + \dots + n_m\vartheta_m)}
$$
of a degree $N \ge 3$ generated by coefficients $a_n(f)$ of functions $f \in S$.  }

\bigskip\noindent
{\bf 5.2. Extremality of Grunsky-Milin functional}. For completeness, we also mention the features of
the general plurisubharmonic functionals on classes of univalent functions on arbitrary quasididsks.

Let again $D$ and $D^*$ be the interior and exterior domains of a bounded (positively oriented) quasicircle $L$. Consider the class $S(D)$ of univalent functions $f$ on $D$ normalized by $f(z) = z + a_2 z^2 + \dots$ near $z = 0$ and the class $\Sigma(D^*)$ of univalent functions $F$ on $D^*$ with hydrodynamical normalization $F(z) = b_0 + b_1 z^{-1} + \dots$ at the infinite point $z = \iy$.
Denote by  $S_k(D)$ and $\Sigma_k(D^*)$ their subclasses formed by functions admitting $k$-quasiconformal extensions across the common boundary $L$ to the whole plane $\hC = \C \cup \{\iy\}$. The unions
$$
S_Q(D) = \bigcup_{k < 1} S_k(D), \quad \Sigma_Q(D^*) =  \bigcup_{k < 1} \Sigma_k(D^*)
$$
are dense in $S(D)$ and $\Sigma(D^*)$ in the topology of locally inform convergence on $D$ and $D^*$,
respectively. The extensions can be additionally normalized, respectively, by $f^\mu(0) = 0$ and $F^\mu(\iy) = \iy$. Both classes are compact with respect to the indicated convergence.

Consider on either from these classes a bounded holomorphic (continuous and Gateaux $\C$-differentiable) functional $J(f)$, which means that for any $f$ and small $t \in \C$,
$$
J(f + t h) = J(f) + t J_f^\prime(h) + O(t^2), \quad t \to 0,
$$
in the topology of uniform convergence on compact sets in $D$ or $D^*$.
Here $J_f^\prime(h)$ is a $\C$-linear functional which is lifted to the strong (Fr\'{e}chet) derivative of $J$ in the norms  of both spaces $L_\iy(D^*)$ and $\B(D)$ (accordingly, in $L_\iy(D)$ and $\B(D^*)$.
Passing if needed to functional $J(f)/M$, one can assume that $M = \max_{S(D)} |J(f)| = 1$.

Varying $f$, one gets the corresponding functional derivative
 \be\label{16}
\psi_0(z) = J_{\id}^\prime(g(\id, z)),
\end{equation}
where
$$
g(w, \z) = \fc{w^2}{w - \z} \quad \text{for} \ \  w \in S(D) \quad \text{and} \ \  g(w, \z) = \fc{1}{w - \z} \quad \text{for} \ \ w \in \Sigma(D^*)
$$
is the kernel of variation on $S_Q(D)$. Any such functional is represented by a complex Borel measure on $\C$ and extends thereby to all holomorphic  functions on $D^*$ (cf. \cite{Sch}). We assume that this derivative is meromorphic on $\C$
and has in the domain $D$ only a finite number of the simple poles (hence $\psi_0$ is integrable over $D$).

For example, one can take the general distortion functionals
 \be\label{17}
J(f) = J(f(z_1), f'(z_1), \dots \ , f^{(\a_1)}(z_1); \dots; f(z_p), f'(z_p), \dots \ , f^{(\a_p)}(z_p))
\end{equation}
with $\grad \wh J(\mathbf 0) \ne 0$, where $\wh J(\mu) = J(f^\mu); \ z_1, \dots \ , z_p$ are distinct fixed points in $D^*$ with assigned orders $\a_1, \dots, \a_p$, and $J$ is a holomorphic function of its arguments.

In this case, $\psi_0$ is a rational function
$$
\wh J_{\id}^\prime(g(\id, z)) = \sum\limits_{j=1}^p \sum\limits_{s=0}^{\a_j-1} \ \fc{\partial
\wh J(\mathbf 0)}{\partial \om_{j,s}} \fc{d^s}{d\z^s} g(w,\z)|_{w=z,\z=z_s},
$$
where $\om_{j,s} = f^{(s)}(z_j)$. Recall that the corresponding Beltrami coefficients run over the unit  ball $\Belt(\D)_1 \subset L_\iy(\D)$.

We associate with the derivative (6) the Teichm\"{u}ller disk
$$
\D(\mu_0) = \{t \mu_0: \ |t| < 1, \ \ \mu_0 = |\psi_0|/\psi_0\}.
$$

The following general theorem solves explicitly the indicated extremal problem on $\Sigma_k(D^*)$ for any such functional and any $k < 1$ (and thereby also for the classes $S_k(D)$, passing to functions $f(z) = 1/F(1/z)$ univalent in domain $D^{-1} = \{z \in \C: \ 1/z \in D\}$ with the same properties).

\bigskip\noindent
{\bf Theorem 6}. {\it If all zeros of the functional derivative $\psi_0(z) = J_{\id}^\prime(g(\id, z))$
on $D^*$ of a given functional (17) are of even order and $\psi_0(0) \ne 0$, then for any $k, \ 0 < k < 1$,
the maximal value of this functional on the any class $\Sigma_k(D^*)$ is attained at the point ${k|\psi_0|/\psi_0}$ on the disk $\D(\mu_0)$, and
$$
\max_{\Sigma_k(D^*)} |J(F^\mu)| = |J(F^{k|\psi_0|/\psi_0})| = k = \a(F^\mu),
$$
with $\a(F)$ given by
$$
\a(F) = \sup \ \Big\{ \Big\vert \iint_D \mu_0(z) \psi(z) dx dy \Big\vert: \ \psi \in A_1^2(D), \ \|\psi\|_{A_1} = 1 \Big\} \quad (z = x + iy).
$$
Here $A_1(D)$ denotes the subspace in $L_1(D)$ formed by integrable holomorphic functions
(quadratic differentials) on $\D$, and $A_1^2(\D)$ is its subset consisting
of $\psi$ with zeros even order in $D$, i.e., of the squares of holomorphic functions.

In addition, the indicated maximal value of $|J(F^\mu)|$ is attained only when this functional coincides  on the disk $\D(\mu_0)$ with the Grunsky-Milin functional
$$
J_{\x^0}(f^\mu) =  \sum\limits_{m,n=1}^\iy \beta_{m,n}(F^\mu) x_m^0 x_n^0
$$
whose defining point $\x^0 = (x_n^0) \in S(l^2)$ is determined by $\mu_0$.
}

\bigskip
Recall that the Grunsky-Milin coefficients $\beta_{m n}$ of functions $F \in \Sigma(D^*)$ are determined from expansion $$
\log \fc{F(z) - F(\z)}{z - \z} = - \sum\limits_{m, n = 1}^\iy
\fc{\beta_{m n}}{\sqrt{m n} \ \chi(z)^m \ \chi(\z)^n},
$$
where $\chi$ denotes a conformal map of $D^*$ onto the disk $\D^*$ so that $\chi(\iy) = \iy, \ \chi^\prime(\iy) > 0$, and that the inequality
$$
\Big\vert \sum\limits_{m,n = 1}^{\iy} \ \beta_{mn} \ x_m x_n \Big\vert \le 1
$$
with any sequence ${\mathbf x} = (x_n)$ from the unit sphere $S(l^2) = \{\x \in l^2: \ \|\x\|= 1\}$
is necessary and sufficient for univalent extension of a function
$f(z) = z + \const + O(1/z)$ near $z = \iy$ to $D^*$
(cf. \cite{Gr}, \cite{Kr7}, \cite{Mi}, \cite{Po1}). The quantity
$$
\vk_{D^*}(f) = \sup \Big\{ \Big\vert \sum\limits_{m,n = 1}^{\iy} \ \beta_{mn} \ x_m x_n \Big\vert : \
{\mathbf x} = (x_n) \in S(l^2)\Big\}
$$
is called the {Grunsky norm} of $f$  and is majorated by the {\bf Teichm\"{u}ller norm} $k(f)$
equal to the smallest dilatation $k(f^{\mu_0})$ among quasiconformal extensions of function $f$ onto $D$;
moreover, on the open dense subset of $\Sigma_Q(D^*)$ the strict inequality $\vk_{D^*}(f) < k(f)$ is valid.

For any quasidisk $D^*$, the set of $f$ with $\vk_{D^*}(f) = k(f)$ is sparse, but such functions play a crucial role in many applications.

In the case of the canonical disk $\D^*$, Theorem 6 provides additionally the values of other basic analytic and geometric quasiinvariants of the extremal curves $L_0 = F^{\mu_0}(|z| = 1)$ (its quasireflection coefficient, Fredholm eigenvalue intrinsically connected with the Grunsky and Teichm\"{u}ller norms of the associated Riemann mapping functions). For details see \cite{Kr12}
and the last section.

\bigskip\bigskip
\centerline{\bf 6. STRENGTHENED VERSION OF THEOREM 1}

\bigskip\noindent
{\bf 6.1}.
Theorem 1 can be essentially strengthened taking more specific families $\mathcal X$ of holomorphic functions on the unit disk considered in Section {\bf 3}.

Assume again that $\mathcal X$ is a collection  of normalized univalent solutions $w(z)$ of equations $S_w = \vp$ with
$$
\vp(z) = \chi(\x) = \sum\limits_0^\iy c_n z^n, \ \quad |z| < 1,
$$
and let, in addition to properties $(a), \ (b)$, this family contains together with each its
function $\vp(z)$ also all functions $\vp_n(z) = \vp(z^n), \ n = 1, 2, \dots$, and the function
$\vp_0 = S_{w_0}$ corresponding to  $w_0(z)$ with  $a_2(w_0) = \max_{\mathcal X} |a_2(w)|$ satisfies
\be\label{18}
\sup_{n \ge 2} |c_n| \le |c_1| \quad (c_1 \ne 0).
\end{equation}

For such families, we have the following general theorem.

\bigskip\noindent
{\bf Theorem 7}. {\it For any family of univalent functions $\mathcal X =\{w\}$, satisfying the above assumptions, the coefficients of their Schwarzians
$$
\vp(z)= S_w(z) = c_0 + c_1 z + \dots + c_n z^n + \dots
$$
are estimated by
 \be\label{19}
|c_n| \le |c_1^0|.
\end{equation}
The equality is attained on function $\vp_0(z^n) = c_0^0 + c_1^0 z^n  + \dots$ arising from
the Schwarzian $\vp_0(z) = S_{w_0}(z)$ of univalent function $w_0$, which maximizes the second
coefficient $a_2$ on $\mathcal X$. }

\bigskip
Applications of this theorem rely on verifying (18), which is complicated.

The assumption (18) has a sense only for subclasses $\mathcal X \ne S$ and imply the bounds of type (19).

Such estimate trivially holds for collection $\mathcal X = B_{r_0}$ considered in Example 1.
The extremal function maximizing $|c_1|$ is $w_0(z)$ whose Schwarzian equals $\sqrt{2/\pi} r_0 z$;
for this family,
$$
\max_{B_{r_0}} \ |c_1| = r_0.
$$

\bigskip\noindent
{\bf 6.2}.
In a similar fashion, one can take the appropriate ball from the weighted space $\mathcal X(\sigma)$
of holomorphic functions $\vp$ in $\D$ with norm
$$
\|\vp\|_{\mathcal X(\sigma)} = \sup_\D (|z| - 1)^\sigma |\vp(z)|
$$
with $0 <\sigma < 2$. The corresponding univalent functions $w$ defined by such $\vp$ are asymptotically conformal on the boundary.

If $0 < \sigma < 1$, then the spaces $\mathcal X(\sigma)$ are embedded into the Bergman spaces $A_p(\D)$,
and for $p = 2 m, \ m > 1$, one obtains explicitly, using Parseval's equality, the subcollections of
functions obeying (18).

The estimate (19) also is valid, for example, for collection of $\vp$ from the unit ball $B_1(H^p)$
of the Hardy space $H^p, \ p > 1$. In this case, it simply follows from Brown's result estimating
$c_1$ for nonvanishing $H^p$ functions (see Lemma 7 below).

\bigskip\bigskip
\centerline{\bf 7. ESTIMATING THE COEFFICIENTS OF SCHWARZIANS}.
\centerline{\bf APLLICATIONS TO NONVANISHING HARDY FUNCTIONS}
\centerline{\bf AND OTHER HOLOMORPHIC FUNCTIONS}

\bigskip
Theorems 1, 2, 7 have many other remarkable applications based on estimating the coefficients of
the Scwarzians.

\bigskip\noindent
{\bf 7.1. Estimating the coefficients of Schwarzians}.
The second application concerns sharp estimation of the coefficients of
Schwarzian derivatives
$$
S_f(z) = \sum\limits_0^\iy \a_n z^n \quad (|z| < 1)
$$
on $S$. This problem also has been investigated by many authors.

For $f_{\tau, \theta}(z) = e^{- i \theta} f(e^{i \tau} z), \ f \in S$, we have
$S_{f_{\tau, \theta}}(z) = e^{2i \theta} S_f(e^{i \tau} z)$. Noting also that the Schwarzian
of the Koebe function is given by
$$
S_{\kp_0}(z) = - \fc{6}{(1 - z^2)^2} = - 6 \ \sum\limits_0^\iy (n + 1) z^{2n},
$$
(and is even), one immediately obtains from Theorem 2 that the even coefficients $\a_{2 n}$ of all Schwarzians $S_f$ on $S$ are sharply estimated by
$$
|\a_{2 n}(S_f)| \le |\a_{2 n}(\kp_0)| = 6(n + 1),
$$
with equality only for $f = \kp_\theta$. Together with the Moebius invariance of $S_f$ and its higher derivatives, the estimate (6) provides the bound of distortion at an arbitrary point of $\D$ (cf., e.g., \cite{Ah}, \cite{Ber}, \cite{H}, \cite{Sci}, \cite{T}).

\bigskip\noindent
{\bf 7.2. The Krzyz and Hummel-Scheinberg-Zalcman conjectures}.
The mean-value inequality producing the relation (8) gives rise to completely different applications of estimates for Schwarzians.

There were two eminent conjectures for nonvanishing  holomorphic functions $f(z)$ on the unit disk $\D$  from the Hardy spaces $H^p$ with $1 < p \le \iy$. Recall that the norm in these spaces is defined by
$$
\|f\|_p = \sup_{r<1} \Bigl(\fc{1}{2 \pi} \ \int\limits_0^{2 \pi} |f(r
e ^{i \theta})|^p d \theta \Bigr)^{1/p}
$$
for $p < \iy$ and $\|f\|_\iy$ coincides with $L_\iy$-norm of this function.

The {\it Krzyz conjecture} \cite{Kz} of 1968 states that Taylor's coefficients of nonvanishing holomorphic functions $f(z) = \sum\limits_0^\iy c_n z^n \in H^\iy$ with $\|f\|_\iy \le 1$
are sharply estimated by $|c_n| \le 2/e$, with equality only for the function
$\kappa_\iy(z^n)$, where
$$
\kappa_\iy(z) = \exp \ \Bigl(\frac{z -1}{z + 1}
\Bigr) = \fc{1}{e} + \fc{2}{e} z - \fc{2}{3e} z^3 + ... \
$$
(and its compositions with rotations about the origin).

Its deep generalization to spaces $H^p$ with $p > 1$ is given by the {\it Hummel-Scheinberg-Zalcman  conjecture} posed in 1977, which states for all nonvanishing $f \in H^p$  with $\|f\|_p \le 1$ the sharp estimate
 \be\label{20}
|c_n| \le (2/e)^{1 - 1/p} ,
\end{equation}
and this bound is realized only by the functions $\epsilon_2
\kappa_{n,p}(\epsilon_1 z)$, where $|\epsilon_1| = |\epsilon_2| = 1$ and
\be\label{21}
\kappa_{n,p}(z) = \Bigl[\frac{(1 + z^n)^2}{2}\Bigr]^{1/p} \ \Bigl[\exp
\frac{z^n - 1}{z^n + 1}\Bigr]^{1-1/p}.
\end{equation}

Both conjectures have been investigated by many authors established only for the initial coefficients
$c_n$  (see, e.g., \cite{Br}, \cite{BK}, \cite{BW}, \cite{Su},\cite{Ho}, \cite{LS}, \cite{PS}, \cite{Sa}, \cite{Sz}, \cite{Su}, \cite{Ta} and the references cited there).
and have been proven only recently by the author in \cite{Kr9}, \cite{Kr10}.

The result of \cite{Kr9}, \cite{Kr10} states:

\bigskip\noindent
{\bf Theorem 8}. {\it The estimate (20) is valid for all spaces $H^p$ with $p \ge 2$; that is, the coefficients of any
nonvanishing function $f \in \ H^p, \ p \ge 2$, with $\|f\|_p \le 1$ satisfy
$|c_n| \le (2/e)^{1 - 1/p}$ for any $n > 1$;  the equality in (1) is realized only on the function $f(z) = \kappa_{n, p}(z)$ and its compositions with pre and post rotations about the origin. }

\bigskip
On going to the limit as $m \to \infty$ one obtains, as a consequence of Theorem 1, that {\it the coefficients of all nonvanishing functions $f(z) = \sum_0^\iy c_n z^n \in H^\iy$ on the unit disk,  with $\|f\|_\iy \le 1$ satisfy the inequality
 \be\label{22}
|c_n| \le \inf \{|c_n(f)|: \ f \in \bigcap\limits_{m \ge 1} \  B_1^0(H^{2m})\}
= \inf \{|c_n(f)|: \ f \in \bigcap\limits_{p \ge 1} \  B_1^0(H^p)\} = 2/e
\end{equation}
for all $n > 1$; here $B_1^0(H^p)$ denotes the collection of nonvanishing
functions from the unit ball in $H^p$.}

This estimate $\max \ |c_n| = 2/e$ for $f \in B_1^0(H^\iy)$ is sharp being realized by the function $\kappa_\iy(z^n)$ and its compositions with rotations. Note that the function $\kappa_\iy(z)$
is the universal holomorphic covering map of the punctured disk $\D \setminus \{0\}$ by $\D$.

This proves the Krzyz conjecture.
Some generalizations of this conjecture also can be proved in this way, however there remain here some
interesting open questions discussed below in subsection {\bf 10.4}.

\bigskip\noindent
{\bf 7.3}.
In fact the estimate (22) is valid on much broader set of functions, because there exist unbounded nonvanishing holomorphic functions $f(z)$ on the unit disk which belong to all spaces $H^p, \ p > 0$.
For example, one can take the functions
$$
f_a(z) = \log \frac{a}{1 - z} \quad \text{with} \ \ |a - 1| \ge 2
$$
(it is a slight modification of a known function which belongs to all $H^p$ with $p > 0$, see \cite{Pr}).

\bigskip\bigskip
\centerline{\bf 8. OUTLINE OF PROOFS OF THEOREMS 7 AND 8}

\bigskip\noindent
{\bf 8.1}. Theorem 8 is a special case covered by Theorems 1 and 7 and, as was mentioned above,
an important step in the proof is to verify the validity of the estimate (18). This is possible
under assumptions of Theorem 8, thus outline here the proof this theorem. The proof of Theorem 7
follows the same lines (having (18)). We deduce the assertion of Theorem 8 as a consequence of
several lemmas.

\bigskip\noindent
$(a)$ \ Denote the unit ball of $H^p$ by $B_1(H^p)$ and its subset of nonvanishing functions by
$B_1^0(H^p)$. It will be convenient to regard the free coefficients $c_0(f)$ also as elements of $B_1^0(H^p)$, which are constant on the disk $\D$. Let
$$
\wh B_1^0(H^p) = B_1^0(H^p) \cup \{f_0\},
$$
where $f_0(z) \equiv 0$. The corresponding sets for the disks $\D_r =\{|z| < r\}$ will be denoted by $B_1^0(H^p(\D_r))$ and $\wh B_1^0(H^p(\D_r))$.

The following lemma concerns the topological features of sets of
nonvanishing holomorphic functions.
Consider the subsets $\mathcal B_r $ of $B_1^0(H^p)$ defined by
$$
\mathcal B_r = \{f \in B_1(H^p(\D_r)): \ f(z) \ne 0 \ \
\text{on the disk} \ \ \D_r = \{|z| < r\}\}, \quad 1 < r < \iy;
$$
then $\mathcal B_{r^\prime} \Subset \mathcal B_r$ if $r^\prime > r$.
Put
$$
B_{*}(H^p) = \bigcup_{r>1} \mathcal B_r
$$
with topology of the inductive limit. All functions from $B_{*}(H^p)$ are zero free in $\D$.

Any point $f_0 \in B_{*}(H^p)$ belongs to all sets $\mathcal B_r$ with
$r \ge r_0$.
Consider the intersection of $B_{*}(H^p)$ with the balls
$\{f \in H^p(B_{r_0}): \ \|f - f_0\|_p < \epsilon\}$ and denote their
connected components containing $f_0$ by $U(f_0, \epsilon)$.

\bigskip\noindent
{\bf Lemma 4}. {\it Each point $f \in B_{*}(H^p)$ has a
neighborhood $U(f, \epsilon)$ in $B_{*}(H^p)$ filled by the functions
which are zero free in the disk $\D$.
Take the maximal neighborhoods $U(f, \epsilon)$ with such property.
Then their union
$$
\mathcal U^p = \bigcup_{f\in B_{*}(H^p)} U(f, \epsilon)
$$
is an open path-wise connective set, hence a domain, in the space $B_{*}(H^p)$.       }

\bigskip
Combining the estimate
$$
\fc{1}{\pi} \ \iint\limits_\D |f(z)|^p dx dy  =
\fc{1}{\pi} \ \int\limits_0^1 \Bigl(\int\limits_0^{2 \pi} |f(r e
^{i \theta})|^p d \theta\Bigr) r dr \le \fc{1}{2} \|f\|_{H^p}^p.
$$
with $\|f\|_\B \le \|f\|_{A_1(\D)}$ (for $f \in A_1$),
one concludes that all functions $f \in H^p$ belong to the space $\B$, and hence these functions (more precisely, the corresponding quadratic differentials
$f(z) dz^2$) can be regarded as the Schwarzian derivatives of locally univalent functions in $\D$.
\footnote{The $H^p$ functions are moved to their equivalence classes so that $f$ and $f_1$ are equivalent if $f_1(z) = \epsilon_1 f(\epsilon_2 z)$ for some constants $\epsilon_1, \ \epsilon_2$ with modulus $1$.
Such equivalence preserves $H^p$ norm and moduli of coefficients. }

\bigskip\noindent
$(b)$ \ The following two lemmas are variations of Lemma 3 and ensure the existence of quasiconformal deformations of $H^p$ functions preserving $H^p$-norm.

First assume that $p = 2 m, \ m \ge 1$ and consider an arbitrary function
$f(z) \in H^{2m} \cap L_\iy(\D)$, with
$$
\sup_\D |f(z)| = M  > \|f\|_{2m}.
$$
Let $E$ be a ring domain bounded by a closed curve $L \subset \D$ containing inside the origin and by the unit circle $S^1 = \partial \D$. The degenerated cases $E = \D \setminus \{0\}$ and $E = S^1$ correspond to the Bergman space $B^p$ and to the Hardy space $H^p$.

Let, in addition,
$$
{\bold d}^0 = (0, 1, 0, ... \ , 0) =: (d_k^0) \in \mathbb R^{n+1},
$$
and $|\mathbf x|$ denote the Euclidean norm in $\mathbb R^l$.

\bigskip\noindent
{\bf Lemma 5}. {\it For any holomorphic function $f(z) =
\sum\limits_{k=j}^\infty c_k^0 z^k \in L_{2m}(E) \cap L_\infty(E)$
(with $c_j^0 \ne 0, \ 0 \le j < n$ and $m \in \mathbb N$), which is not a polynomial of degree
$n_1 \le n$, there exists a positive number $\ve_0$ such that for every point
$$
{\mathbf d}^\prime = (d'_{j+1}, \dots, \ , d'_n) \in \C^{n-j}
$$
and every $a \in \mathbb R$ satisfying the inequalities
$$
|{\mathbf  d}^\prime| \le \ve, \ \ |a| \le \ve,
\ \ \ve < \ve_0,
$$
there exists a quasiconformal automorphism $h$ of the complex plane $\hC$, which is conformal in the disk
$$
D_0 = \{w : \ |w - c_0^0| < \sup_\D |f_0(z)| + |c_0^0| +1\}
$$
(hence also outside of $E$) and satisfies the conditions:

(i) $h^{(k)}(c_0^0) = k! d_k = k! (d_k^0 + d_k^\prime), \  k = j +
1, \dots, n$ (i.e., $d_1 = 1 + d_1^\prime$ and
$d_k = d_k^\prime$ for $k \ge 2)$;

(ii) $\|h \circ f \|_{2m}^{2m} =  \|f \|_{2m}^{2m} + a$. }

\bigskip
Note that quasiconformal deformations of such type preserve the $L_p$-norm of holomorphic functions, but generically increase their $L_\iy$-norm.

The proof of this lemma shows that all its assumptions are essential; the arguments do not extend to arbitrary $p \ge 2$ and unbounded holomorphic $L^p$ functions.

As for arbitrary $p > 1$, the indicated arguments can be applied to {\bf nonvanishing} functions from the corresponding spaces $H^p$, and this leads to the following weakened extension of Lemma 5.

\bigskip\noindent
{\bf Lemma 6}. {\it For every bounded nonvanishing function $f(z) = c_0 + c_1 z + c_2 z^2 + \dots \in H^p, \ p > 1$, with $\|f\|_p < \iy$, which is not a polynomial of degree at most $n$, there exists $\ve_0 > 0$ such that for any point $\mathbf d \in \mathbb C^{n+1}$ with $|\mathbf d| \le \ve < \ve_0$, there is a bounded nonvanishing function $f^*(z) = c_0^* + c_1^* z + c_2^*  z^2 + \dots \in H^p$ satisfying $c_j^* = c_j + d_j$ for all $j = 0, 1, \dots, n$ and }
$$
\|f^*\|_p = \|f\|_p + O(\ve).
$$

The point is that for any nonvanishing function $f(z) = \sum\limits_0^\iy c_n z^n \in H^p, p > 1$,
the function
$$
f_{p/2}(z) := f(z)^{p/2} = e^{(p/2) \log f(z)},
$$
with a fixed branch of the logarithmic function, also is single valued, holomorphic and zero free in
the unit disk $\D$; we take the principal branch.

We also use the following important Brown's result \cite{Br} quoted above, which provides the needed estimate of the first coefficient

\bigskip\noindent
{\bf Lemma 7}. \cite{Br} {\it For any $f(z) = c_0 + c_1 z + c_2 z^2 + \dots \in B_1^0(H^p)$, we have
$$
|c_1| \le (2/e)^{1 - 1/p},
$$
with equality only for the rotations of function $\kappa_{1,p}(z)$ given by (21
). }

$(c)$ \ Similar to Theorem 2, each functionals $J_n(f) = c_n$ are lifted to the universal Teichm\"{u}ller space $\T$ and its cover $\T_1$ producing on $\Fib(\T)$ amaximal subharmanic function $u_n(t)$.
This is established similat to Theorem 2, with some technical differences (one has to use the approximation of functions $f \in H^p$ and by polynomials and, (to apply Lemmas 5 and 6) by bounded functions.

To estimate the seconde coefficient $c_2$, one must maximize on the set $B_1^0(H^{2m})$ the functional
$$
I_2(f) = \max \ (|J_2(f)|, |J_2(f_2)|) = \max \ (|c_2(f)|, |c_2(f_2)|),
$$
where $f_9z) = f(z^2)$.
Since the correspondence $f(z) \mapsto f_2(z)$ is linear, it is holomorphic in $H^p$ norm;
thus this functional also is plurisubharmonic.

One can repeat for this functional the above construction, lifting both functionals
$J_2(f)$ and $J_2(f_2)$ onto the space $\T_1$ similar to above, and obtain in the same
way the corresponding nonconstant radial subharmonic function on the disk
$\D_{2|a_2^0|} = \{w: \ |w| < 2 |a_2^0|\}$.
Again, this function is logarithmically convex, hence monotone increasing, and attains its maximal value at $|t| = 2 |a_2^0|$.

Now observe that by Parseval's equality for the boundary functions
$$
f(e^{i \theta}) = \lim_{r\to 1} f(r e^{i \theta}), \quad f \in H^{2m},
$$
we have
$$
1 \ge \frac{1}{2 \pi} \int\limits_{-\pi}^\pi |f(e^{i \theta})|^2 d \theta
= \sum_1^\infty |c_n|^2.
$$
Applying it to the function
$$
f(z) = \kappa_{1,m}(z) = \sum\limits_0^\iy c_n^0 z^n
$$
and noting that by (21),
$$
|c_1^0|^2 = (2/e)^{2(1-1/(2m))} = 0.5041...^{1-1/(2m)} > 0.5041... \ ,
$$
one obtains that for this function,
\be\label{23}
\sum\limits_2^\infty |c_n^0|^2 < 0.5 < |c_1^0|^2,
\end{equation}
which provides the desired estimate (20) for $n = 2$.

Taking subsequently the functions $f_3(z) = f(z^3), \ f_4(z) = f(z^4), \dots$
and the corresponding functionals
$$
\begin{aligned}
I_3(f) &= \max \ (|J_3(f)|, |J_3(f_3)|) = \max(|c_3(f)|, |c_3(f_3)|),  \\
I_4(f) &= \max \ (|J_4(f)|, |J_4(f_4)|) = \max(|c_4(f)|, |c_4(f_4)|), \dots \ ,
\end{aligned}
$$
one obtains by the same arguments that the estimate (20) is valid for all $n \in \mathbb N$.

\bigskip\noindent
{\bf 8.2. Remark on the case $1 < p < 2$}.
All arguments in the proof of Theorem 7, excluding the Parseval equality, applied in the last step, work for any $p > 1$.
In fact, this equality was applied only to the function $\kappa_{1,p}$ maximizing
$|c_1|$ and was used for estimation $|c_n|$ by comparison of the initial non-free coefficient of functions $\kappa_{1,p}(z^m), \ 1 \le m \le n$.

The explicit representation (21) of $\kappa_{1,p}$ shows that this function also is bounded on the unit disk for $p$ satisfying  $1 < p < 2$; hence it belongs to $H^2$. However, for all such $p$, we have
$$
\|\kappa_{1,p}\|_{H^2} > \|\kappa_{1,p}\|_{H^p};
$$
thus the needed relation (23) and its applications in the last step $(b)$ fail. Therefore, the functionals $J_n$ and $I_n$ cannot be compared on this way.

\bigskip\noindent
{\bf 8.3}. The above remark also shows that estimating the coefficients $c_n$ with $n > 1$ of functions
$S_w = \chi(\x)$ arising in Theorem 7 essentially depends on the inner features of this family.

\bigskip\bigskip
\centerline{\bf 9. REMARKS ON EXTREMAL FUNCTIONS IN GENERAL BANACH}
\centerline{\bf SPACES}

\bigskip
One of the interesting extensions of the Hummel-Scheinberg-Zalcman problem (also still unsolved) is
to estimate the Taylor coefficients of nonvanishing holomorphic maps $f(z) = c_0 + c_1 z + \dots$
of the unit disk $\D$ into other complex Banach spaces $X$.
Denote by $B(X)$ the unit ball of $X$.

The following remarks concern two special cases of this general problem and provide its restricted solutions.

First, one can show following \cite{Kr9} that the collection of all nonvanishing functions $\vp \in \B$ forms a Banach submanifold, which allows one to apply the arguments exposed above and estimate
non-zero coefficients of such Schwarzians.

The second case concerns the Bergman's space $A_2$, for which we show that the features of extremal functions can be essentially different from the Hardy spaces. Recall that the norm of $A_2$ is  $\|f\| = (\fc{1}{\pi} \iint_\D |f(z)|^2 dx dy)^{1/2}$.

The collection $B_0(A_2)$ of nonvanishing holomorphic functions $f(z)$  mapping the
disk $\D$ into the closed ball $\ov{B(A_2)}$ (i.e., with $f(z) \ne 0$ on $\D$ and $\|f\| \le 1$) is compact
in the weak topology of the locally uniform convergence in $\D$. So any holomorphic coefficient functional
 \be\label{24}
J(f) = J(c_{m_1}, \dots, c_{m_s}) \quad \text{with} \ \ 1 \le m_1 < m_2 < \dots < m_s = N < \iy
\end{equation}
has an extremal $f_0$ on which $|J(f)|$ attains its maximum on $B_0(A_2)$.

While the extremal functions of many problems in Hardy spaces are bounded,
Proposition 2 implies, for example, that {\it any function $f_0 \in B_0(A_2)$ maximizing the functional (24) must be unbounded on the disk $\D$} (compare with the extremal problems
for nonvanishing Bergman functions investigated e.g. in \cite{ABKS}, \cite{BK}).

This difference is caused by the fact mentioned after the proof of Proposition 2: quasiconformal deformations
created by Propositions 1 and 2 preserve the norm in $A^p$, while the norm of the Hardy spaces can be increased.

Indeed, it follows from Proposition 2 that any extremal $f_0$ of $J(f)$ on $B_0(A_2)$  must be unbounded on $\D$, unless $f_0$ is a zero-free polynomial
 \be\label{25}
p_N(z) = c_0 + c_1 z + \dots + c_N z^N
\end{equation}
(with $c_0 \ne 0$); otherwise, one can vary the coefficients $c_k$ and obtain by this lemma an admissible function $f_{*} \in B_0(A_2)$ with $|J(f_{*})| > |J(f_0)|$.

It remains to establish that the polynomials (25) with $\|p_N\|_{A_2} \le 1$ cannot be extremal for $J(f)$.
We pick a sufficiently small $\ve > 0$ and consider
the polynomial
$$
p_\ve(z) = -\ve c_0 + \ve z^{N+1},
$$
for which
$$
\max_{\mathbb S^1} |p_\ve(z)| < \max_{\mathbb S^1} |p_N(z)|.
$$
Then the Rouch\'{e} theorem yields that the polynomial
$$
P_{N+1,\ve}(z) = p_N(z) + p_\ve(z) = (1 - \ve) c_0 + c_1 z + \dots + c_N z^N + \ve z^{N+1}
$$
also must be, together with $p_N$, zero-free on $\D$. Its norm is estimated by
$$
\|p_{N+1,\ve}\|_{A_2}^2 = (1 - \ve)^2 |c_0|^2 + \fc{|c_1|^2}{2} + \dots
+ \fc{|c_N|^2}{N + 1} + \fc{\ve^2}{N + 2} = \|p_N\|_{A^2}^2 - 2 \ve + O(\ve^2) <   \|p_N\|_{A^2}^2 = 1,
$$
which implies that $P_{N+1,\ve}$ is an admissible function, with
 \be\label{26}
|J(p_{N+1,\ve})| = |J(p_N)| = \max \{|J(f)|: \ f \in B_0(A^2)\}.
\end{equation}
But this contradicts to Proposition 2, because this proposition allows one to construct the variations of $p_{N+1,\ve}$, which preserve its $A^2$-norm and increase $|J(p_{N+1,\ve})|$, disturbing (26).

Note that this also is a special case of the general Theorem 7.

\bigskip\bigskip
\centerline{\bf 10. UNIVALENCE AND POLYNOMIAL APPROXIMATION}

\bigskip\noindent
{\bf 10.1}.
An interesting not investigated aspect of univalence concerns the speed of approximation of univalent functions $f$ by univalent polynomials, especially in the strong $\B$-norm.
The rigidity of conformality causes that set of corresponding Schwarzians $S_f$ is rather thin in the space $\T$.

There are deep results in complex geometric function theory that the Riemann conformal mapping function of the unit disk $\D$ onto any simply connected hyperbolic domain $D \subset \hC$ can be approximated by {\bf univalent} polynomials, and generically this approximation is uniform on compact subsets of this disk; see, e.g., \cite{MG}. A simple proof of this fact in a strengthened form is given in \cite{Kr11}; it involves
the integral approximation of holomorphic functions and their quasiconformal extension.

First we consider the univalent functions $f$ in the disk $\D$ which are {\bf asymptotically conformal} on the boundary, i.e., map the unit circle $\mathbb S^1$ onto the asymptotically conformal Jordan curves $L$, which means that for any pair of points $a, b \in L$, we have
$$
\max\limits_{z \in L(a,b)} \frac{|a - z| + |z - b|}{|a - b|} \to 1 \quad \text{as} \quad |a - b| \to 0,
$$
where the point $z$ lies on $L$ between $a$ and $b$.

Such curves are quasicircles without corners and can be rather pathological
(see, e.g., \cite{Po2}, p. 249). All $C^1$-smooth curves are asymptotically conformal.

Any univalent function on $\D$ is approximated locally uniformly  by asymptotically conformal functions
by taking the homotopy $f_t(z) = \fc{1}{t} f(t z)$ with $t$ close to $1$.

The main results here are given by the following two theorems from \cite{Kr11}.

\bigskip\noindent
{\bf Theorem 9}. {\it Every univalent function $f(z)$ on $\D$ with asymptotically conformal boundary values  is approximated by univalent polynomials $p_n$ on $\ov \D$ so that $\|S_{p_n} - S_f\|_\B \to 0$ as $n \to \iy$ and $p_n$ admit $k_n$-quasiconformal extensions to $\hC$ with dilatations
$k_n \to k$; in addition, the Grunsky norms $\vk(p_n) \to \vk(f)$. }

\bigskip
Such approximation connected with quasiconformalily holds for the generic univalent functions in
a weakened form.

\bigskip\noindent
{\bf Theorem 10}. {\it For any univalent function $f$ in the disk $\D$ admitting
quasiconformal extension across the unit circle, there exists a sequence of univalent polynomials $p_n$ on the close disk $\ov \D$ convergent to $f$ uniformly on $\ov \D$ and such that their
dilatations $k(p_n) \nearrow k(f)$.     }

These results have important applications, for example, to qualitative characterization of conformal mapping functions. We illustrate this on the Fredholm eigenvalues of curves.

The Fredholm eigenvalues $\rho_n$ of an oriented smooth closed Jordan curve $L  \subset \hC$ are the eigenvalues of its double-layer potential, or equivalently, of the integral equation
 \be\label{27}
u(z) +  \frac{\rho}{\pi} \int\limits_L \ u(\zeta) \frac{\partial}{\partial
n_\zeta} \log \frac{1}{|\zeta - z|} ds_\zeta = h(z),
\end{equation}
where $n_\zeta$ denotes the outer normal and $ds_\zeta$ is the length element at $\zeta \in L$.
These values have crucial applications in solving many problems in various fields of mathematics.

The least positive eigenvalue $\rho_L = \rho_1$ is naturally connected with conformal and quasiconformal maps related to $L$ and can be defined for any oriented closed Jordan curve $L$ by
$$
\fc{1}{\rho_L} = \sup \ \fc{|\mathcal D_G (u) - \mathcal D_{G^*} (u)|} {\mathcal D_G (u) + \mathcal D_{G^*} (u)},
$$
where $G$ and $G^*$ are, respectively, the interior and exterior of $L; \ \mathcal D$ denotes the Dirichlet integral,
and the supremum is taken over all functions $u$ continuous on $\hC$ and harmonic on $G \cup G^*$.

Due to the K\"{u}hnau-Schiffer theorem, the value $\rho_L$ is reciprocal to the Grunsky norm $\vk(f_L^*)$ (see \cite{Ku2}, \cite{Sc2}).

The equation (27) is valid for polynomial curves $p_n (\mathbb S^1)$ and one can apply the known approximal numerical methods for solving (24) (see, e.g. \cite{Ga}). In the limit, this gives the estimates for more general univalent functions.
Note that generically the quasicircles $f(\mathbb S^1)$ are fractal curves.

The Schwarzians of polynomials  $p(z) = \sum\limits_0^n a_j z^j$ with $a_1 \ne 0$
are rational functions on $\hC$ of the form
 \be\label{28}
S_p(z) = \sum\limits_j \fc{c_j}{(z - z_j)^2} \quad \text{with} \ \ \sum\limits_j |c_j| > 0,
\end{equation}
where $z_j$ are the finite critical points of $p$ (i.e., the zeros of its derivative $p^\prime(z)$, which is a polynomial of degree $n - 1$). In the case of univalent $p$ on $\D$, these points are placed outside of $\D$. For such polynomials, we have

\bigskip\noindent
{\bf Theorem 11}. {\it Let a polynomial $p(z)$ be univalent on the disk $\D$ and all finite critical points $z_1, \dots, z_{n-1}$ are placed  on the boundary circle $\mathbb S^1$, and let $\|S_p\|_\B < 2$. Then
 \be\label{29}
k(p) = \vk(p) = q_{p(\mathbb S^1)} = 1/\rho_{p(\mathbb S^1)} = \max_j |c_j|,
\end{equation}
where $c_j$ arise from (28).
}

\bigskip\noindent
Here $q_L$ denotes the reflection coefficient of quasicircle $L$.

\bigskip\noindent
{\bf 10.2}. If the polynomials approximate a function $f \in S_Q$ on the boundary circle $\mathbb S^1$
enough fast, then this function extends holomorphically to a broader disk $\D_R$ with $R > 1$.
This follows from the classical Bernstein-Walsh-Siciak theorem (see, e.g., \cite{Wa}, \cite{Kr6}), which yields that, letting
$$
e_m (h, K) = \inf \big\{ \max\limits_{z \in K} |h(z) - p(z)|: \ p \in \mathcal P_m \big\},
$$
where ${\mathcal P}_m$ is the space of polynomials whose degrees do not exceed $m, \ m = 1, 2, ...$,
a continuous function $h$ on a compact $K \subset C$ extends holomorphically to the region
$$
D_R = \{ z \in {\mathbb C}^n : \ g_K(z) < \log R \}, \quad R > 1,
$$
determined by the Green function $g_K$ of this compact with pole at the infinite point and $R$ satisfying
$$
\limsup\limits_{m \to \iy} e_m^{1/m} (h, K) = 1/R.
$$

\bigskip\noindent
{\bf 10.3. Example 3}. Consider again the weighted space $\mathcal X(\sigma)$ introduced
in {\bf 6.2}.
It is embedded into $\B$, and the univalent solutions $w(z)$ of the equations
$S_w = \vp$ with $\vp \in \mathcal X(\sigma)$ are asymptotically conformal in $\D$, since
$$
\lim\limits_{|z| \to 1-} (1 - |z|^2)^2 S_w (z) = 0.
$$
For all such $w(z)$, we have the relations (29), and their Schwarzians $\vp$ are approximated in $\B$-norm
by fractions (28).

Theorems 1 and 7 yield some estimates for such $w(z)$.

\bigskip\noindent
{\bf 10.4}.
An interesting natural extension of Krzyz's conjecture is to estimate the coefficients of holomorphic maps from the unit disk into its ring subdomain bounded by a closed curve $L \subset \D$ containing inside the origin and the unit circle $\mathbb S^1$.
Even for the rings $G_r = \{r < |z| < 1\}, \ r > 0$, the estimates more precise that (22), are unknown.

Another question related to this conjecture is to find the bound for coefficients of $\vp$ from the unit ball in the space $\mathcal X(\sigma)$ with $0 < \sigma < 1$.

\bigskip\bigskip
\centerline{\bf 11. BIUNIVALENT FUNCTIONS}

\bigskip\noindent
{\bf 11.1}. Let us mention here an interesting class of univalent functions for which we have a complete distortion theory.
It consists of so-called {\bf biunivalent} functions, that is of $w = f(z) \in S$ whose inverse functions
$z = f^{-1}(w)$ also are univalent on the disk $\D$. Biunivalence causes rather strong rigidity, and for such functions the boundary of the image domain $f(\D)$ touches the unit circle outside.

There is an intrinsic connection of biunivalent functions with special functions, solutions of complex differential equations, with the so-called $q$-calculus, etc., wher biunivalence naturally arises. From these points of view, biunivalent functions have been
and remain be intensively investigated by many authors, which considered and defined new special
subclasses of such functions depending on different parameters;
see, e.g., \cite{AlH}, \cite{BT}, \cite{BrO}, \cite{FA}, \cite{Le}, \cite{Ne}, \cite{PK}, \cite{SMG},
\cite{SrB} and the references cited there. These investigations resulted mainly in the estimates of the
initial Taylor coefficients $a_2$ and $a_3$ and of some their combinations.

\bigskip\noindent
{\bf 11.2. General distortion theorem}. Put
 \be\label{30}
|a_2(\mathcal B)| = \sup \{|a_2(f)|: \ f \in \mathcal B\},
\end{equation}
and consider the set of rotations
$$
\mathcal R_\mathcal B  = \{f_{0,\tau,\theta}(z) = e^{- i \theta} f_0(e^{i \tau} z)\},
$$
where $f_0$ is one of the maximizing functions for $a_2$ on $\mathcal B$ (their existence follows from compactness of $\mathcal B$ in topology of locally uniform convergence in $\D$).

E. Netanyahu established that $|a_2(\mathcal B)| = 4/3$ \cite{Ne}.
We shall call any maximizing function for (30) the Netanyahu function.

\bigskip
The general theorem is rather surprising (in view of the rigidity of biunivalence) though its assertion
is similar to the above distortion coefficient theorems.

\bigskip\noindent
{\bf Theorem 12}. {\it Any rotationally invariant polynomial functional (4), whose zero set
$\mathcal Z_J = \{f \in \mathcal B: \ J(f) = 0\}$
is separated from the rotation set $\mathcal R_\mathcal B$, is maximized only by the Netanyahu functions $f_{0,\tau,\theta}$.  }

\bigskip
One has to establish first that the class $\mathcal B$ admits the needed properties (openness of
its image in the universal space $\T$ and variational stability) underlying the proofs of the previous  distortion theorems. The second condition is fulfilled rather trivially. Openness is a consequence of
the following result from \cite{Kr13} and from complex homogeneity of the universal space $\T$.

\bigskip\noindent
{\bf Theorem 13}. \cite{Kr13} {\it Any holomorphic function $f(z)$ on the disk $\D$, whose Schwarzian
$S_f$ lies in the ball $B_\T(\mathbf 0, \kappa)$ of radius $\kappa \le 1/4$ admits $\kappa$-quasiconformal extension onto $\hC$ and is of the form
$$
f(z) = \gamma \circ f_0(z),
$$
where $f_0$ is biunivalent on $\D$ and $\gamma$ is a Moebius transformation of $\hC$.
The upper bound $1/4$ is sharp (cannot be increased).  }

\bigskip
{\bf 11.3. Biunivalent functions with $k$-quasiconformal extension}.
Another remarkable fact related to biunivalent functions is that (in contract to collections of all normalized univalent functions) the subclasses $\mathcal B_k$ of $\mathcal B$ consisting of functions with $k$-quasiconformal extensions with a fixed  $k < 1$ also obey both  properties (conditions) $(a), \ (b)$. This allows one  to estimate sharply the rotationally invariant coefficient functionals on any class $\mathcal B_k$.

\bigskip\noindent
{\bf Theorem 14}. {\it Let $w = f_k(z)$ be a maximizing function for
$|a_2(\mathcal B_k)| = \sup_{f \in \mathcal B_k} |a_2(f)|$.
Then any rotationally invariant polynomial functional (4) on $\mathcal B_k$,
whose zero set $\mathcal Z_J$ is separated from the rotation set $f_{k,\tau,\theta}$, is maximized only by these rotations.
}

\bigskip
In particular, this theorem provides the sharp estimates of coefficients $a_n$ of $f \in \mathcal B_k$
and coefficients $\a_n$ for all nonzero coefficients of the Netanyahu function.

The point is that we need  for such functions the quasiconformal variations which do not increase the maximal dilatation $k$. Such variations were constructed in the book  \cite{Kr1}
in terms of inverse maps $z = f^{-1}(w)$. Biunivalence yields that both functions $f$ and $f^{-1}$ simultaneously belong to $\mathcal B_k$, which allows us to apply such variation to the original
functions $f$ (cf. \cite{Kr12}). It can be also shown that these variations with appropriately small dilatations preserve biunivalence.

The details will be presented elsewhere.

\bigskip\bigskip
\centerline{\bf 12. SOME OPEN PROBLEMS}

\bigskip\noindent
{\bf 1}. Find the extent to which the restrictions of Theorem 2 can be weakened.

\bigskip\noindent
{\bf 2}. Find the sufficient conditions, different from (18), for arbitrary families of univalent functions, which allow to estimate the coefficients of their Schwarzians.

\bigskip\noindent
{\bf 3}. Nothing is known on speed of approximation of univalent functions by univalent polynomials.
It would be very interesting to obtain some quantitative results.

\bigskip\noindent
{\bf 4}. The above theory is intrinsically connected with the classes of univalent functions
having standard quasiconformal extensions (satisfying the classical Beltrami equation).

Find the extent to which the results of such type are valid for univalent functions admitting more general  qusiconformal extensions, for example, the extensions which belong to the Gutlyanskii-Ryazanov classes of extremal solutions of quasilinear and more general Beltrami equations (considered in \cite{GR}).

\bigskip
{\small\it{ \leftline{Department of Mathematics, Bar-Ilan
University, 5290002 Ramat-Gan, Israel} \leftline{and
Department of Mathematics, University of Virginia,  Charlottesville, VA 22904-4137, USA}}

\end{document}